\documentclass[11pt]{article}

\usepackage{mathrsfs}
\usepackage{graphicx}
\usepackage{amsmath}
\usepackage{amsfonts}
\usepackage{amssymb}

\def\sqr#1#2{{\vcenter{\vbox{\hrule height.#2pt
              \hbox{\vrule width.#2pt height#1pt \kern#1pt \vrule width.#2pt}
              \hrule height.#2pt}}}}
\def\signed #1{{\unskip\nobreak\hfil\penalty50
              \hskip2em\hbox{}\nobreak\hfil#1
              \parfillskip=0pt \finalhyphendemerits=0 \par}}
\def\endpf{\signed {$\sqr69$}}

\def\ns{\noalign{\medskip}}

\newtheorem{lemma}{Lemma}[section]
\newtheorem{theorem}{Theorem}[section]

\newtheorem{proposition}{Proposition}[section]

\newtheorem{definition}{Definition}[section]

\newtheorem{remark}{Remark}[section]

\makeatletter
   
   \@addtoreset{equation}{section}
\makeatother
\date{}

\begin{document}

\title{A nonhomogeneous boundary value problem for the Kuramoto-Sivashinsky equation in a quarter plane\thanks {This work was partially supported by the
National Natural Science Foundation of China (No. 11301425,  No.
11571244 and No. 11231007) and the State Scholarship Fund (No.
201306985014) from China Scholarship Council.
\medskip}}

\author{Jing Li\thanks{School of Economics and Mathematics, Southwestern
University of Finance and Economics, Chengdu 611130,
China.{\small\it E-mail:} {\small\tt jingli@swufe.edu.cn}.\medskip}
\quad Bing-Yu Zhang\thanks{Department of Mathematical Sciences,
University of Cincinnati, Cincinnati, 45221, USA. {\small\it
E-mail:} {\small\tt bzhang@math.uc.edu}.\medskip} \quad and \quad
Zhixiong Zhang\thanks{School of Mathematics, Sichuan University,
Chengdu 610064, China. {\small\it E-mail:} {\small\tt
zxzhang@amss.ac.cn}.}}

\date{}

\maketitle

\begin{abstract}
{\noindent We study the initial boundary value problem for
one-dimensional Kuramoto-Sivashinsky equation with nonhomogeneous
boundary conditions. Through the analysis of the boundary integral
operator, and applying the known results on the Cauchy problem, we
obtain both the local well-posedness and the global well-posedness
for the nonhomogeneous initial boundary value problem. It is shown
that the Kuramoto-Sivashinsky equation is well-posed in Sobolev
space $C([0,T]; H^s (\mathbb{R}^+)) \bigcap L^2(0,T;
H^{s+2}(\mathbb{R}^+))$ for $s>-2$. }
\end{abstract}

%

\bigskip


\bigskip

\section{Introduction}

More than thirty years ago, Kuramoto and Sivashinsky derived the
Kuramoto-Sivashinsky equation independently by different methods.
Sivashinsky \cite{Sivashinsky1,Sivashinsky2} obtained an asymptotic
nonlinear equation which describes the evolution of a disturbed
plane flame front. Kuramoto discussed the turbulent phase waves in
\cite{Kuramoto2} and the behavior of a curved wavefront in
\cite{Kuramoto3}, and got the analogous equation in both cases.

There are many works addressing the Kuramoto-Sivashinsky equation,
see \cite{Aimar, Biagioni, Cerpa2, Cousin, Kaikina, Nicolaenko1,
Pilod, Tadmor} for the well-posedness problem, see \cite{Collet,
Demirkaya, Giacomelli, Ilyashenko, Nicolaenko2} for the long-time
behavior, see \cite{Armaou, Cerpa1, Cerpa2, Glowinski} for the
control problem, and the rich references cited therein.

As mentioned above, the well-posedness problem of the
Kuramoto-Sivashinsky equation in one space dimension has been
studied extensively. Most of the references on well-posedness
problems cited above are addressed to either periodic boundary
conditions or homogeneous boundary conditions.
B. Nicolaenko and B. Scheurer (\cite{Nicolaenko1}) exploited
standard energy methods to discuss the initial value problem with
periodic boundary conditions. In \cite{Tadmor}, E. Tadmor proved the
existence and stability for the pure Cauchy problem and the periodic
problem. A global solution was obtained by patching together short
time solutions, which does not require a  higher regularity priori
estimate as energy methods do.
H. A. Biagioni and T. Gramchev (\cite{Biagioni}) discussed the
initial value problem problem for multidimensional
Kuramoto-Sivashinsky type equations on $\mathbb{R}^n$ and on the
torus $\mathbb{T}^n$. The initial data could belong to $H^r_p$ with
$r$ negative. \cite{Biagioni} also obtained global well-posedness in
$L^2$  when the nonlinearities are conservative.
In \cite{Cousin}, A. T. Cousin and N. A. Larkin considered the
initial boundary value problem with homogeneous boundary values in a
non-cylindrical domain. Using the Faedo-Galerkin method, they proved
the existence and uniqueness of global weak, strong and smooth
solutions.
In \cite{Kaikina}, E. I. Kaikina studied the global existence and
large-time behavior of solutions to the homogeneous boundary value
problem for subcritical Kuramoto-Sivashinsky-type equation on a
half-line with initial data in $L^{\infty} \cap L^{1,2}$ being
sufficiently small, where $L^{1,2}$ is some weighted Lebesgue space.
D. Pilod (\cite{Pilod}) was interested in the Cauchy problem for the
dispersive Kuramoto-Velarde equation with initial datum in
$H^s(\mathbb{R})$ ($s>-1$). In Section 2 of \cite{Cerpa2}, because
of the control framework they considered, E. Cerpa and A. Mercado
got the well-posedness results for the boundary data $u(0,t),
u_x(0,t)\in L^2(0,T)$, $u(1,t)= u_x(1,t)=0$, and initial datum
$u(x,0)\in H^{-2}(0,1)$.

In this paper, we study the nonhomogeneous boundary value problem
for the following dispersive Kuramoto-Sivashinsky equation
(sometimes, we call it Korteweg-de Vries-Kuramoto-Sivashinsky
equation) in a quarter plane:
\begin{eqnarray}\label{1.1.1*}
\left\{\begin{array}{ll} \displaystyle u_{t} + u_{xxxx} + \delta
u_{xxx}+ u_{xx} + uu_x=0,  &  (x,t)\in \mathbb{R}^+ \times (0,T),
\\ \ns
u(x,0)=\phi (x),   &  x\in \mathbb{R}^+,
\\ \ns
u(0,t)=h_1(t), \quad  u_x(0,t)=h_2(t),  \qquad & t\in (0,T),
\end{array}\right.
\end{eqnarray}
where $T>0$ and $\delta\in \mathbb{R}$.
For initial datum  $\phi \in H^s(\mathbb{R}^+)$, one may wonder what
is the lowest regularity of the boundary conditions to guarantee
solution $u \in C ([0, T]; H^s(\mathbb{R}^+))\cap L^2  ( 0,T ;
H^{s+2} (\mathbb{R}^+) )$. One may also be interested in what is the
smallest value of $s$ that guarantee the well-possedness of equation
(\ref{1.1.1*}). Theorem \ref{Thm 1.2*} (see below) gives answers to
these two questions in some sense. The optimal regularity for the
boundary conditions are
 $(h_1, h_2) \in
 H^{\frac{s}{4}+\frac{3}{8}} (0,T) \times H^{\frac{s}{4}
+\frac{1}{8}} (0,T) $.  It is shown that when $s>-2$, the solution
exists for an arbitrary long time. In \cite{Pilod}, D. Pilod proved
that the Kuramoto-Velarde equation
$$ \displaystyle \partial_t w +
\delta \partial_x^3 w + \mu( \partial_x^4 w  + \partial_x^2 w  )  =
\alpha ( \partial _x w )^2  +  \gamma  w  \partial_x^2 w $$
was globally well-posed in $H^s (\mathbb{R})$ for $s
> -1$ and this result was sharp in the sense that the flow map of the equation fails
to be $C^2$ in $H^s(\mathbb{R})$ when $s<-1$. Taking $\mu=1$,
$\alpha = - \frac{1}{2}$ and $\gamma=0$, and setting $u=w_t$, one
can obtain the sharp well-posedness result for $s
> -2$ concerning the Cauchy problem:
$$u_{t} + u_{xxxx} + \delta u_{xxx} + u_{xx} + uu_x=0.$$
We guess that $s>-2$ is also
sharp in our nonhomogeneous boundary value problem case. In order to
prove the sharpness, one need to construct a counterexample which
implies that the wellposedness will fail when $s< -2$. This has not
been done yet in this paper.

Our main task in this paper is to yield the well-posedness result
for $-2 < s < 0$. Also, we obtain the well-posedness result for $s
\ge 0$. These results are presented in a unified way in Theorem
\ref{Thm 1.2*} below. When $s>\frac{1}{2}$, one need the following
definition of $s$-compatibility conditions.

\begin{definition}\label{Def 1.1*}
Let $T>0$ and $s\geq 0$ be given. $(\phi, h_1, h_2) \in  H^s
(\mathbb{R}^+) \times H^{\frac{s}{4}+\frac{3}{8}} (0,T) \times
H^{\frac{s}{4} +\frac{1}{8}} (0,T)$ is said to be $s$-compatible if
\begin{enumerate}
\item[i)] When $s - 4[\frac{s}{4}] \leq \frac{1}{2}$,
then
$$\begin{array}{ll}  \phi_k (0) = h_1 ^{(k)} (0),  \quad  \phi_k^{\prime} (0) = h_2 ^{(k)} (0),  & \quad  k=0,1,\cdots, [\frac{s}{4}] - 1.
\end{array}$$

\item[ii)] When $\frac{1}{2} < s - 4[\frac{s}{4}] \leq
\frac{3}{2}$, then
$$\begin{array}{ll}   \phi_k (0) = h_1 ^{(k)} (0), \quad  \phi_k^{\prime} (0) = h_2 ^{(k)} (0),  & \quad  k=0,1,\cdots, [\frac{s}{4}] - 1;
\\ \ns
   \phi_k (0) = h_1 ^{(k)} (0),  & \quad   k= [\frac{s}{4}] .
\end{array}$$

\item[iii)] When $ s - 4 [\frac{s}{4}] >
\frac{3}{2}$, then
$$ \begin{array}{ll}  \phi_k (0) = h_1 ^{(k)} (0),  \quad  \phi_k^{\prime} (0) = h_2 ^{(k)} (0),  & \quad   k=0,1,\cdots, [\frac{s}{4}] .
\end{array}$$

\end{enumerate}
Here,
$$\left\{\begin{array}{ll} \displaystyle  \phi_0 (x) = \phi(x),
\\ \ns
 \displaystyle \phi_k (x) = - \phi_{k-1}^{\prime\prime\prime\prime} (x) - \delta \phi_{k-1}^{\prime\prime\prime} (x) - \phi_{k-1}^{\prime\prime} (x)  - \sum _{j=0} ^{k-1}  C_{k-1}^j \phi_j(x) \phi_{k-j-1} ^{\prime}(x)
 ,\quad  k=1,2, \cdots
\end{array}\right.$$
Recall that $C_{k-1}^j  =\frac{(k-1)!} {j! (k-j-1)!}$ and $[\cdot]$
is the floor function.
\end{definition}

Our main result in this article is stated as follows.

\begin{theorem}\label{Thm 1.2*}
Let  $T>0$ and $\varepsilon>0$.
\begin{enumerate}
\item[i)] If $s\geq 0$ and $(\phi, h_1, h_2) \in  H^s (\mathbb{R}^+)
\times H^{\frac{s}{4}+\frac{3}{8}} (0,T) \times H^{\frac{s}{4}
+\frac{1}{8}} (0,T)$ is $s$-compatible, then equation (\ref{1.1.1*})
admits a unique solution $u\in C ([0, T]; H^s( \mathbb{R}^+ ))\cap
L^2  ( 0,T ; H^{s+2} ( \mathbb{R}^+ ) )$ with
$ u_x \in C (  [0,+\infty); H^{\frac{s}{4} + \frac{1}{8}} (0,T) )$.
Moreover, the corresponding solution map from the space of initial
and boundary data to the solution space is  continuous.

\item[ii)] If   $-2  < s<
0$,
  $\phi \in H^s (\mathbb{R}^+)$,   $(h_1, h_2) \in
 H^{\frac{s}{4}+\frac{3}{8}} (0,T) \times H^{\frac{s}{4}
+\frac{1}{8}} (0,T) $   and $( t^{\frac{|s|}{4} + \varepsilon} h_1,
t^{\frac{|s|}{4} + \varepsilon} h_2 ) \in H^{ \frac{3}{8}} (0,T)
\times H^{ \frac{1}{8}} (0,T)$, then equation (\ref{1.1.1*}) admits
a unique solution $u\in C ([0, T]; H^s( \mathbb{R}^+ )) $.
Moreover, the corresponding solution map from the space of initial
and boundary data to the solution space is   continuous.
\end{enumerate}

\end{theorem}

To prove the above theorem, firstly, by the method developed in
\cite{Zhang1, Zhang2, Zhang3} which is addressed to the Korteweg-de
Vries equation, we discuss the nonhomogeneous boundary value problem
of the associated linear equation. The key technical part is to
transfer the homogeneous initial boundary value problem
 to the Cauchy problem, through a careful
study of the boundary integral operator $W_{bdr}(t)$ (see ``Proof of
estimate (\ref{2.2.2*})"). Here, $v = W_{bdr} (t) (h_1 (t), h_2(t))$
is the solution of
\begin{eqnarray*}
\left\{\begin{array}{ll} \displaystyle v_{t} + v_{xxxx} +  \delta
v_{xxx} = 0,  & (x,t)\in \mathbb{R}^+ \times (0,T),
\\ \ns
v(x,0)=0,   &  x\in \mathbb{R}^+,
\\ \ns
v(0,t)=h_1(t), \quad  v_x(0,t)=h_2(t),  \qquad & t\in (0,T).
\end{array}\right.
\end{eqnarray*}
And $v$ can be explicitly expressed by the Laplace transform (see
Remark \ref{Rem 2.1}).

Secondly, thanks to the bilinear estimates and the fixed point
theory, we prove the existence of a local solution for the
Kuramoto-Sivashinsky equation (\ref{1.1.1*}).

Finally, combining the local well-posedness results and a priori
estimate, we can patch local solutions together to a global
solution. Hence, the Kuramoto-Sivashinsky equation
 (\ref{1.1.1*}) is globally well-posed.

%

The rest of this paper is organized as follows. In Section 2, we
consider the smoothing properties of the linear Cauchy problem.
Section 3 is devoted to discussing the nonhomogeneous boundary value
problem  of the associated linear equation. Then we consider the
well-posedness  of  Kuramoto-Sivashinsky equation (\ref{1.1.1*}) for
$s\geq 0$ and $-2< s <0$ in Section 4 and Section 5, respectively.

%


\section{The linear problem in a whole line}\label{Sec 2}

In this section, we consider the smoothing properties of the linear
Cauchy problem.

\subsection{Results on the Cauchy problem}\label{Sec 2.1}


\begin{proposition}\label{Pro 6.1}
Let $ s \in \mathbb{R}$, $T>0$ and $\phi \in  H^s (\mathbb{R})$.
Then equation
\begin{eqnarray}\label{2.2.7}
\left\{\begin{array}{ll} \displaystyle p_{t} + p_{xxxx} + \delta
p_{xxx} =0, \qquad & (x,t)\in \mathbb{R}\times (0,T),
\\ \ns
p(x,0)= \phi (x),  &   x\in \mathbb{R},
\end{array}\right.
\end{eqnarray}
admits a unique solution $p  \in C ( [0,T] ; H^s (\mathbb{R}) ) \cap
L^2 ( 0,T ; H^{s+2} (\mathbb{R})  )$  with
$\partial^k_x p \in C ( \mathbb{R} ; H^{\frac{s-k}{4} + \frac{3}{8}}
(0,T) )$ ($k=0,1$).
Moreover, there exist a constant $C> 0$ such that for any $\phi \in
H^s(\mathbb{R})$, it holds
\begin{equation}\label{2.2.8}
\sup_{t\in [0,T]} \| p (\cdot, t) \| _{H^s (\mathbb{R})}  \leq C \|
\phi \| _{H ^s (\mathbb{R})},      \qquad   \| p\| _{L^2( 0,T;
H^{s+2} (\mathbb{R}) )}    \leq  C \| \phi \| _{H ^s (\mathbb{R})},
\end{equation}
and
\begin{equation}\label{2.2.10}
\sup _{x\in \mathbb{R}}  \| \partial ^k _x  p(x, \cdot) \|
_{H^{\frac{s-k}{4} + \frac{3}{8}} (0,T)}
 \leq C  \| \phi \| _{ H
^s(\mathbb{R}) }.
\end{equation}
Furthermore, if $s<0$ and $0<T\leq 1$, it holds
\begin{equation}\label{2.2.11}
\sup _{t \in [0,T] }   \|  t^{\frac{|s|} {4} }    p( \cdot, t ) \|
_{L^2 (\mathbb{R})} \leq C \| \phi \| _{ H ^s(\mathbb{R}) }, \qquad
\| t^{\frac{|s|} {4} }  p\| _{L^2( 0,T; H^{2} (\mathbb{R}) )}
  \leq C  \| \phi \| _{ H ^s(\mathbb{R}) }
\end{equation}
and
\begin{equation}\label{2.2.12}
\sup _{x\in \mathbb{R}}  \| t^{\frac{|s|} {4} } \partial ^k _x p(x,
\cdot) \| _{H^{-\frac{k}{4} + \frac{3}{8}} (0,T)}
 \leq C  \| \phi \| _{ H
^s(\mathbb{R}) }.
\end{equation}

\end{proposition}

\begin{proposition}\label{Pro 6.2}
Let $s\in \mathbb{R}$, $T>0$ and $f(x,t) \in L^1 (0,T; H ^s (
\mathbb{R} ) )$. Then equation
\begin{eqnarray}\label{2.4.5}
\left\{\begin{array}{ll} \displaystyle p_{t} + p_{xxxx} + \delta
p_{xxx} = f(x,t),
 \qquad  & (x,t)\in
\mathbb{R}\times (0,T),
\\ \ns
p(x,0)= 0,  &  x\in \mathbb{R},
\end{array}\right.
\end{eqnarray}
admits a unique solution $p  \in C ( [0,T] ; H^s (\mathbb{R}) ) \cap
L^2 ( 0,T ; H^{s+2} (\mathbb{R})  )$  with
$\partial^k_x p \in C ( \mathbb{R} ; H^{\frac{s-k}{4} + \frac{3}{8}}
(0,T) )$ ($k=0,1$).
Moreover, there exist a constant $C> 0$ such that for any $f(x,t)
\in L^1 (0,T; H ^s ( \mathbb{R} ) )$, it holds
\begin{equation}\label{2.4.7}
\sup_{t\in [0,T]} \| p (\cdot, t) \| _{H^s (\mathbb{R})}  \leq  C
\|f \|_{L^1 (0,T; H^s(\mathbb{R}))},    \quad     \| p\| _{L^2( 0,T;
H^{s+2} (\mathbb{R}) )}    \leq  C  \|f \|_{L^1 (0,T;
H^s(\mathbb{R}))},
\end{equation}
and
\begin{equation}\label{2.4.9}
\sup _{x\in \mathbb{R}}  \| \partial ^k _x  p(x, \cdot) \|
_{H^{\frac{s-k}{4} + \frac{3}{8}} (0,T)} \leq C  \|f \|_{L^1 (0,T;
H^s(\mathbb{R}))}.
\end{equation}
Furthermore, if $s<0$, $0<T\leq 1$ and $  t^{\frac{|s|} {4} }
f(x,t) \in L^1 (0,T; L^2 ( \mathbb{R} ) )$, it holds
\begin{equation*}\label{2.4.11}
\sup _{t \in [0,T] }  \|  t^{\frac{|s|} {4} }   p( \cdot, t ) \|
_{L^2 (\mathbb{R})}
+ \| t^{\frac{|s|} {4} }  p\| _{L^2( 0,T; H^{2} (\mathbb{R}) )}
 \leq
C  \left( \|f \|_{L^1(0,T; H^s(\mathbb{R}))} + \|t^{\frac{|s|} {4} }
f \|_{L^1 (0,T; L^2(\mathbb{R}))}   \right)
\end{equation*}
and
\begin{equation}\label{2.4.12}
\sup _{x\in \mathbb{R}}  \| t^{\frac{|s|} {4} } \partial ^k _x p(x,
\cdot) \| _{H^{-\frac{k}{4} + \frac{3}{8}} (0,T)}
 \leq
C \left( \|f \|_{L^1(0,T; H^s(\mathbb{R}))} + \|t^{\frac{|s|} {4} }
f \|_{L^1 (0,T; L^2(\mathbb{R}))}   \right).
\end{equation}

\end{proposition}

The proof of Proposition \ref{Pro 6.1} - Proposition \ref{Pro 6.2}
will be given in Subsection \ref{Sec 2.2} - Subsection \ref{Sec
2.3}, respectively.


\subsection{Proof of Proposition \ref{Pro 6.1}}\label{Sec 2.2}

To prove Proposition \ref{Pro 6.1}, we need the following lemma,
which follows from minor modifications of Lemma 2.5 in
\cite{Zhang2}.

\begin{lemma}\label{Lma 2.1.1}
Let $\gamma (\rho)$ be a continuous complex-valued function defined
on $(0,+\infty)$ satisfying the following two conditions:
\begin{enumerate}
\item[i)]There exist $\delta>0$ and $b>0$ such that $\displaystyle \sup_{0<\rho <\delta }  \frac{|{\text Re} \gamma(\rho)|}  {\rho}   \geq b$;

\item[ii)]There exists a complex number $\alpha+i \beta $ such that $\displaystyle \lim _{\rho\rightarrow +\infty} \frac{|\gamma(\rho)|}  {\rho}   =\alpha +i \beta$.
\end{enumerate}
Then there exists a constant $C>0$ such that for all $f\in L^2(0,
+\infty)$,
$$ \left \| \int_0^{+\infty} e ^{\gamma (\rho) t} f(\rho) d \rho \right \| _{L^2(0, T)}  \leq C \left(  \|e^{T {\text Re} \gamma(\cdot)} f(\cdot)\| _{L^2(\mathbb{R}^+)}  + \|f(\cdot)\|_{L^2(\mathbb{R}^+)}  \right).$$
\end{lemma}

\noindent  {\bf Proof of Proposition \ref{Pro 6.1}.} The proof is
divided into several steps.

\smallskip

{\it Step 1.} To prove estimate (\ref{2.2.10}), we only need to show
that for any $r\in \mathbb{ R}$,
\begin{equation}\label{2.2.9}
\sup _{x\in \mathbb{R}} \| \partial ^k _x  p(x, \cdot) \| _{H^r
(0,T)} \leq C \| \phi \| _{H ^{4r +k-\frac{3}{2}} (\mathbb{R})}.
\end{equation}

Firstly, we claim that estimate (\ref{2.2.9}) holds for any
non-negative integer $r$.
In fact,  the solution of equation (\ref{2.2.7}) can be expressed as
follows
\begin{equation}\label{2.2.15}
 p(x,t)  =
\frac{1}{2\pi}\int_{- \infty}^{+ \infty}
 e^{i x \xi} e^{ ( i \delta \xi^3 - \xi^4 ) t} \widehat{\phi}(\xi)
 d\xi,
\end{equation}
where ``$\ \widehat{} \ $'' denotes the Fourier transform with
respect to $x$, i.e.
$$
\widehat{ \phi } (\xi)=\int_{-\infty}^{+\infty} e^{-i x \xi} \phi(x)
dx.$$
Formula (\ref{2.2.15}) gives that for $k=0,1$,
\begin{equation}\label{2.2.17}
\begin{array}{ll}
\partial_x^k    p(x,t) &  \displaystyle  =
\frac{1}{2\pi}\int_{0}^{+ \infty}  i^k \xi^{k}
 e^{i x \xi} e^{ (i \delta \xi^3  - \xi^4 ) t} \widehat{\phi}(\xi)
 d\xi     +     \frac{1}{2\pi}\int_{- \infty}^{0}  i^k \xi^{k}
 e^{i x \xi} e^{ ( i \delta \xi^3  - \xi^4 ) t} \widehat{ \phi }(\xi)
 d\xi    \\  \ns
& \displaystyle  :=  \partial_x^k   p^+(x,t)  + \partial_x^k
p^-(x,t).
 \end{array}
\end{equation}
Setting $\zeta = \xi^4$ in (\ref{2.2.17}), one has
\begin{equation}\label{2.2.40}
\begin{array}{ll}
\displaystyle  \partial_x^k p^+ (x, t) = \frac{1}{8 \pi}\int_{0}^{+
\infty} i^k \sqrt[4]{\zeta ^{k-3}}    e^{i x \sqrt[4]{\zeta}} e^{ (i
\delta \sqrt[4]{\zeta^3}  - \zeta) t} \widehat{\phi}
 (\sqrt[4]{\zeta})
 d\zeta  ,
\end{array}
\end{equation}
which yields that for any $x\in\mathbb{R}$,
\begin{equation}\label{2.2.41}
\begin{array}{ll}
\displaystyle  \|  \partial_x^k p^+ (x, \cdot )  \| ^2 _{H^r (0,T)}
\displaystyle  = \sum_{m=0}^r  \|
\partial_t ^m  \partial_x^k   p(x, \cdot )  \| ^2 _{L^2(0,T )}  \\ \ns
\displaystyle  = \sum_{m=0}^r  \left \|
 \frac{1}{8 \pi}\int_{0}^{+ \infty}  i^k (i
\delta \sqrt[4]{\zeta^3}  - \zeta )^m  \sqrt[4]{\zeta ^{k-3}}   e^{i
x \sqrt[4]{\zeta}}   e^{(i \delta \sqrt[4]{\zeta^3}  - \zeta ) t}
 \widehat{\phi}
 (\sqrt[4]{\zeta})
 d\zeta  \right  \| ^2 _{L^2(0,T)}  .
\end{array}
\end{equation}
It follows from (\ref{2.2.41}) and Lemma \ref{Lma 2.1.1} that for
any $x\in\mathbb{R}$
\begin{equation}\label{2.2.42}
\begin{array}{ll}
\displaystyle  \|  \partial_x^k p^+ (x, \cdot )  \| ^2 _{H^r
(0,T)}\\ \ns
 \displaystyle  \leq  C  \sum_{m=0}^r    \int_{0}^{+ \infty}
\left( e^{- \zeta T} +1 \right)^2 \left|   i ^k  (i \delta
\sqrt[4]{\zeta^3}  - \zeta )^m  \sqrt[4]{\zeta ^{k-3}}    e^{i x
\sqrt[4]{\zeta}} \widehat{\phi}
 (\sqrt[4]{\zeta})  \right|^2
 d\zeta \\ \ns
 \displaystyle  \leq  C  \sum_{m=0}^r    \int_{0}^{+ \infty} ( \delta^{2m} \sqrt{\zeta^{3m} } + \zeta
^{2m}) \sqrt{\zeta ^{k-3}}    \left| \widehat{\phi}
 (\sqrt[4]{\zeta})  \right|^2
 d\zeta .
\end{array}
\end{equation}
Setting $\zeta = \xi^4$ in (\ref{2.2.42}), we have
\begin{equation}\label{2.2.47}
\begin{array}{ll}
\displaystyle \sup_{x\in \mathbb{R}} \|  \partial_x^k p^+ (x, \cdot
) \| ^2 _{H^r (0,T)}
& \displaystyle  \leq  C  \sum_{m=0}^r    \int_{0}^{+ \infty} (
\delta^{2m} \xi^{6m} + \xi^{8m} ) \xi ^{2k-3}   \left|
\widehat{\phi}
 (\xi)  \right|^2
 d\xi \\ \ns
& \displaystyle  \leq  C    \int_{0}^{+ \infty} (1+\xi^2) ^{4r
+k-\frac{3}{2}} \left| \widehat{\phi}
 (\xi)  \right|^2
 d\xi .
\end{array}
\end{equation}
Similarly, we can obtain
\begin{equation}\label{2.2.48}
\displaystyle \sup_{x\in \mathbb{R}} \|  \partial_x^k p^- (x, \cdot
) \| ^2 _{H^r (0,T)}
 \displaystyle  \leq  C    \int_{- \infty}^0 (1+\xi^2) ^{4r
+k-\frac{3}{2}} \left| \widehat{\phi}
 (\xi)  \right|^2
 d\xi.
\end{equation}
Combining (\ref{2.2.17}), (\ref{2.2.47}) and (\ref{2.2.48}), we
obtain the desired estimate  (\ref{2.2.9}) for any non-negative
integer $r$.

\smallskip

Secondly, we claim that estimate (\ref{2.2.9}) holds for any
negative integer $r$. Indeed, for any $\varphi(t) \in H_0^{-r}
(0,T)$, integrating by parts, we deduce from (\ref{2.2.40}) that
\begin{equation*}
\begin{array}{ll}
\displaystyle \langle   \partial_x^k p^+  (x,t), \varphi(t) \rangle
_{H^r(0,T) \times H_0^{-r}(0,T) }\\ \ns
  \displaystyle   = \int _0^T \left( \frac{1}{8 \pi}\int_{0}^{+
\infty} i^k   \sqrt[4]{\zeta ^{k-3} }   e^{i x \sqrt[4]{\zeta}} e^{
(i \delta \sqrt[4]{\zeta^3}  - \zeta) t}   \widehat{\phi}
 (\sqrt[4]{\zeta})
 d\zeta    \right)  \varphi(t)  dt\\ \ns
 \displaystyle  =  \int _0^T  (-1)^r  \left(
  \frac{1}{8 \pi}\int_{0}^{+
\infty}  i^k  (i \delta \sqrt[4]{\zeta^3}  - \zeta)^r
\sqrt[4]{\zeta ^{k-3} }    e^{i x \sqrt[4]{\zeta}}  e^{ (i \delta
\sqrt[4]{\zeta^3} - \zeta) t} \widehat{\phi}
 (\sqrt[4]{\zeta})
 d\zeta    \right)  \varphi^{(-r)}(t)
   dt,
\end{array}
\end{equation*}
which implies that for any $x\in \mathbb{R}$,
\begin{equation*}
\displaystyle  \|  \partial_x^k p^+ (x, \cdot )  \| ^2 _{H^r (0,T)}
  \displaystyle  \leq  \left\|   \frac{1}{8 \pi} \int_{0}^{+
\infty}  i^k  (i \delta \sqrt[4]{\zeta^3}  - \zeta)^r \sqrt[4]{\zeta
^{k-3} }    e^{i x \sqrt[4]{\zeta}}  e^{ (i \delta \sqrt[4]{\zeta^3}
- \zeta) t} \widehat{\phi}
 (\sqrt[4]{\zeta})
 d\zeta     \right\| ^2_{L^2 (0,T) } .
\end{equation*}
Similar to the proof from (\ref{2.2.42}) to (\ref{2.2.48}), we can
easily get that  estimate  (\ref{2.2.9}) holds for any negative
integer $r$.

\smallskip

Finally, for non-integer values of $s$, estimate  (\ref{2.2.9}) can
be obtained by standard interpolation theory.

\medskip

{\it Step 2.} We claim that (\ref{2.2.8}) holds. In fact, it follows
from (\ref{2.2.15}) that
\begin{equation*}
\begin{array}{ll}
\displaystyle \sup_{t\in[0,T]} \left\| p(\cdot,t) \right\|
_{H^s(\mathbb{R})}
\displaystyle =  \sup_{t\in[0,T]}  \left( \int_{-\infty}^{+\infty}
(1+\xi^2)^{s} \left| e^{ ( i \delta \xi^3 - \xi^4 ) t}
\widehat{\phi} (\xi) \right|^2
  d\xi  \right)^{\frac{1}{2}}  \\ \ns
\displaystyle \leq  \sup_{t\in[0,T]}     \left\| e^{ - \xi^4  t}
\right\|_{L^{\infty} (\mathbb{R})}   \left\|\phi
  \right\|_{H^s(\mathbb{R})}
\leq  \left\|\phi
  \right\|_{H^s(\mathbb{R})},
\end{array}
\end{equation*}
and
\begin{equation*}\label{2.2.50}
\begin{array}{ll}
\displaystyle \left\| p \right\|^2 _{L^2(0,T; H^{s+2} (\mathbb{R}))}
\displaystyle = \int_0^T   \int_{-\infty}^{+\infty} (1+\xi^2)^{s+2}
\left| e^{ ( i \delta \xi^3 - \xi^4 ) t}  \widehat{\phi} (\xi)
\right|^2
  d\xi   dt \\ \ns
\displaystyle  =  \int_{-\infty}^{+\infty}  \left(  \int_0^T
 (1+\xi^2)^{2} e^{ -2  \xi^4  t}   dt \right) (1+\xi^2)^{s}
\left| \widehat{\phi} (\xi) \right|^2  d\xi  \\ \ns
\displaystyle  \leq   \int_{-\infty}^{+\infty}  \left(  \int_0^T
\left( 2+2\xi^4  e^{- 2\xi^4t }  \right)  dt \right) (1+\xi^2)^{s}
\left| \widehat{\phi} (\xi) \right|^2  d\xi \\ \ns
\leq C \left\|\phi
  \right\|^2 _{H^s(\mathbb{R})}.
\end{array}
\end{equation*}

{\it Step 3.} We claim that (\ref{2.2.11}) and (\ref{2.2.12}) are
true. Firstly, when $t\in[0,T] \subset [0,1]$, it holds $
t^{\frac{|s|} {2} } \leq (1+ t^{\frac{1}{2}} \xi^2)^{|s|}
(1+\xi^2)^{-|s|}$. Then by (\ref{2.2.15}),
\begin{equation*}
\begin{array}{ll}
\displaystyle   \sup_{t\in[0,T]}    \| t^{\frac{|s|} {4} }
p(\cdot,t)  \| _{L^2(\mathbb{R})}
\displaystyle =   \sup_{t\in[0,T]}  \left( \int_{-\infty}^{+\infty}
t^{\frac{|s|}{2}}  \left| e^{ ( i \delta \xi^3 - \xi^4 ) t}
\widehat{\phi} (\xi) \right|^2
  d\xi  \right)^{\frac{1}{2}}  \\ \ns
\displaystyle  \leq   \sup_{t\in[0,T]}   \left(
\int_{-\infty}^{+\infty}   (1+ t^{\frac{1}{2}} \xi^2)^{|s|} e^{ -2
\xi^4  t}   (1+\xi^2)^{s} \left|\widehat{\phi} (\xi) \right|^2
  d\xi  \right)^{\frac{1}{2}} \\ \ns
\displaystyle  \leq 2^{|s|-1}   \sup_{t\in[0,T]}  \left(   \left\|
e^{- \xi^4t }  \right\|_{L^{\infty} (\mathbb{R})}   +   \left \| (t
\xi^4)^{\frac{|s|} {4} } e^{- \xi^4t }  \right\|_{L^{\infty}
(\mathbb{R})} \right) \left\|\phi
  \right\|_{H^s(\mathbb{R})} \\ \ns
\displaystyle  \leq 2^{|s|-1} \left (1+ \frac{|s|} {4} ^{\frac{|s|}
{4} } e^{{\frac{s} {4} }} \right) \left\|\phi
  \right\|_{H^s(\mathbb{R})}
 \leq 2^{|s|} \left\|\phi
  \right\|_{H^s(\mathbb{R})}.
\end{array}
\end{equation*}

Secondly, by (\ref{2.2.7}), $q_1:= tp$ satisfies
\begin{eqnarray}\label{2.2.53}
\left\{\begin{array}{ll} \displaystyle  \partial_t q_{1} +
\partial_x^4 q_{1} + \delta \partial_x^3 q_{1} =p, \qquad & (x,t)\in \mathbb{R}\times (0,T),
\\ \ns
q_1 (x,0)= 0, &  x\in \mathbb{R},
\end{array}\right.
\end{eqnarray}
where  $p \in  L^2 ( 0,T ; H^{s+2} (\mathbb{R}) ) $. Thus,  maximal
regularity property implies that the solution of equation
(\ref{2.2.53}) satisfies
\begin{equation}\label{2.2.56}
\|tp \|_{ L^2 ( 0, T ; H^{s+6} (\mathbb{R}) ) } =   \|q_1 \|_{ L^2 (
0, T ; H^{s+6} (\mathbb{R}) ) }    \leq C  \|p \|_ {L^2 ( 0,T ;
H^{s+2} (\mathbb{R}) )}.
\end{equation}
Similarly, by (\ref{2.2.7}), $q_2:= t^2 p$ satisfies
\begin{eqnarray}\label{2.2.53*}
\left\{\begin{array}{ll} \displaystyle  \partial_t q_{2} +
\partial_x^4 q_{2} + \delta \partial_x^3 q_{2} = 2tp, \qquad & (x,t)\in \mathbb{R}\times (0,T),
\\ \ns
q_2 (x,0)= 0, &  x\in \mathbb{R},
\end{array}\right.
\end{eqnarray}
where  $t p \in  L^2 ( 0,T ; H^{s+6} (\mathbb{R}) ) $. Thus, maximal
regularity property and (\ref{2.2.56}) imply that the solution of
equation (\ref{2.2.53*}) satisfies
\begin{equation*}\label{2.2.56*}
\|t^2 p \|_{ L^2 ( 0, T ; H^{s+10} (\mathbb{R}) ) }  \leq C  \|p \|_
{L^2 ( 0,T ; H^{s+2} (\mathbb{R}) )}.
\end{equation*}
Iterating the argument, we have that for any $k\in \mathbb{Z}^+$,
\begin{equation}\label{2.2.56**}
\|t^k p \|_{ L^2 ( 0, T ; H^{s+2+4k} (\mathbb{R}) ) }  \leq C  \|p
\|_ {L^2 ( 0,T ; H^{s+2} (\mathbb{R}) )}.
\end{equation}
The Young's inequality implies that for any $s<0$, $t\in [0,T]$,
$\xi \in \mathbb{R}$ and for some fixed $k>\frac{|s|}{4}$,
$$ (1+\xi^2) ^{2}   t^{\frac{|s|}{2}}  =  (1+\xi^2) ^{\frac{(s+2)(s+4k)}{4k}} \cdot  (1+\xi^2) ^{ \frac{|s|(s+2+4k)}{4k}}   t^{\frac{|s|}{2}}
\leq  \frac{s+4k}{4k} (1+\xi^2) ^{s+2} +  \frac{|s|}{4k}  (1+\xi^2)
^{s+2+4k} t^{2k}.$$
Therefore, combing (\ref{2.2.8}) and (\ref{2.2.56**}), we have
\begin{equation*}
\begin{array}{ll}
\displaystyle   \| t^{\frac{|s|}{4} } p \|^2 _{L^2(0,T;
H^2(\mathbb{R}))}  =  \int_0^T  \left( \int_{-\infty}^{+\infty}
(1+\xi^2) ^{2}  t^{\frac{|s|}{2}} \left |\widehat{ p } (\xi,t)
\right|^2 d\xi \right)   dt
\\ \ns
\displaystyle   \leq  \frac{s+4k}{4k}   \int_0^T  \left(
\int_{-\infty}^{+\infty} (1+\xi^2) ^{s+2}  \left |\widehat{ p }
(\xi,t) \right|^2 d\xi \right) dt \\ \ns
\displaystyle  \quad +   \frac{|s|}{4k} \int_0^T  \left(
\int_{-\infty}^{+\infty} (1+\xi^2) ^{s+2+4k} t^{2k} \left
|\widehat{p} (\xi,t) \right|^2 d\xi \right) dt
\\ \ns
\displaystyle  \leq   \left \| p \right\|^2 _{ L^2 ( 0, T ; H^{s+2}
(\mathbb{R}) )}
+    \left \| t^k  p \right\|^2 _{ L^2 ( 0, T ; H^{s+2+4k}
(\mathbb{R}) )}
\leq C  \|\phi \|^2_{H^s(\mathbb{R})} .
\end{array}
\end{equation*}

Finally,  by (\ref{2.2.11}), $ t^{\frac{|s|}{4} } p \in C([0,T];
L^2(\mathbb{R}))  \cap   L^2(0,T; H^2(\mathbb{R})) $. Trace Theorem
shows that (\ref{2.2.12}) holds.

This completes the proof of Proposition \ref{Pro 6.1}.\endpf

\begin{remark}\label{Rem 2.2}

Let $ A: {\cal D}(A)\subset L^2(\mathbb{R}^+) \rightarrow
L^2(\mathbb{R}^+)$ be a linear operator defined as follows
\begin{equation*}\label{2.2.0}
A w = -w_{xxxx}, \qquad {\cal D}( A ) = H^4(\mathbb{R}^+)\cap
H_0^2(\mathbb{R}^+).
\end{equation*}
Then $ A $ is the infinitesimal generator of an analytic semigroup
in $L^2(\mathbb{R}^+)$. With the perturbation of a lower order
operator, ${\cal A} w := A w - \delta w_{xxx}$ is the infinitesimal
generator of a $C_0$ semigroup $W_c(t)$ in $L^2(\mathbb{R}^+)$.

Let $\bar{{\cal A}}: H^4(\mathbb{R}) \rightarrow L^2(\mathbb{R})$ be
a linear operator defined as follows
\begin{equation*}\label{2.2.5}
\bar{{\cal A}} w = -w_{xxxx}- \delta w_{xxx}.
\end{equation*}
Then $\bar{\cal A}$ is the infinitesimal generator of a $C_0$
semigroup
$W_{\mathbb{R}}(t)$ in $L^2(\mathbb{R})$. %
Moreover, the solution of equation (\ref{2.2.7}) with initial datum
$\phi$ can be expressed by
\begin{equation}\label{2.2.21}
W_{\mathbb{R}}(t) \phi   =  \frac{1}{2\pi}\int_{- \infty}^{+ \infty}
 e^{i x \xi}  e^{ ( i \delta \xi^3 - \xi^4 ) t}  \widehat{\phi}(\xi)
 d\xi.
\end{equation}

\end{remark}

\subsection{Proof of Proposition \ref{Pro 6.2}}\label{Sec 2.3}

The proof is divided into several steps.

\smallskip

{\it Step 1.} We claim that (\ref{2.4.9}) is true. Indeed, by
semigroup theory (\cite{Pazy}), the solution of equation
(\ref{2.4.5}) can be formally written as follows
\begin{equation}\label{2.4.15}
p(x,t) = \int_0^t W_{\mathbb{R}} (t-\tau) f(x, \tau)  d\tau    =
\int_0^T   \chi (\tau)  W_{\mathbb{R}} (t-\tau) f(x, \tau) d\tau,
\end{equation}
where
$$\chi (\tau)=  \left\{\begin{array}{ll} 1, \quad
\tau \leq t,
\\ \ns
0, \quad \tau>t.
\end{array}\right.$$
Thus, it follows from (\ref{2.4.15}) that for $k=0,1$,
\begin{equation}\label{2.4.17}
\partial ^k _x  p(x, t)   = \int_0^T   \chi
(\tau)
\partial ^k _x (W_{\mathbb{R}} (t-\tau) f(x, \tau) ) d\tau.
\end{equation}
Hence, by (\ref{2.4.17}) and the Minkowski's inequality, we see that
for any $x\in \mathbb{R}$,
\begin{equation}\label{2.4.19}
\begin{array}{ll}
\displaystyle  \| \partial ^k _x  p(x, \cdot) \|_{H^{\frac{s-k}{4} +
\frac{3}{8}} (0,T)}  = \left \| \int_0^T
 \chi (\tau)  \partial ^k _x  (W_{\mathbb{R}}  (t-\tau) f(x, \tau) ) d\tau \right
\|_{H^{\frac{s-k}{4} + \frac{3}{8}} (0,T)}  \\  \ns
\displaystyle  \leq   \int_0^T   \left \|  \chi(\tau) \partial ^k _x
(W_{\mathbb{R}}  (t-\tau) f(x, \tau) ) \right \|_{H^{\frac{s-k}{4} +
\frac{3}{8}} (0,T)}  d\tau  \\  \ns
\displaystyle  =  \int_0^T   \left \| \partial ^k _x (W_{\mathbb{R}}
(t-\tau) f(x, \tau) ) \right \|_{H^{\frac{s-k}{4} + \frac{3}{8}}
(\tau,T)} d\tau.
\end{array}
\end{equation}
Proposition \ref{Pro 6.1} and  (\ref{2.4.19}) show that
\begin{equation*}
\begin{array}{ll}
\displaystyle \sup _{x\in \mathbb{R}} \| \partial ^k _x  p(x, \cdot)
\| _{H^{\frac{s-k}{4}  +  \frac{3}{8}} (0,T)}
\leq   \int_0^T  \sup _{x\in \mathbb{R} }    \left \|
\partial ^k _x (W_{\mathbb{R}} (t-\tau) f(x, \tau) ) \right
\|_{H^{\frac{s-k}{4} + \frac{3}{8}} (\tau,T)}  d\tau   \\  \ns
\displaystyle  \leq C  \int_0^T   \|f(\cdot, \tau)\|_{H
^s(\mathbb{R})} d\tau
= C \| f \| _{L^1(0,T; H ^s(\mathbb{R}))}.
\end{array}
\end{equation*}

{\it Step 2.} We claim that
\begin{equation}\label{2.4.20}
\sup _{t \in [0,T] }   \|  t^{\frac{|s|} {4} }   p( \cdot, t ) \|
_{L^2 (\mathbb{R})}  \leq 2^{ \frac{5|s|}{4} -1 }   \left( \|f
\|_{L^1(0,T; H^s(\mathbb{R}))} + \|t^{\frac{|s|} {4} } f \|_{L^1
(0,T; L^2(\mathbb{R}))}   \right).
\end{equation}
Indeed, by (\ref{2.4.15}) and the Minkowski's inequality, we see
that
\begin{equation}\label{2.4.21}
\begin{array}{ll}
\displaystyle   \sup_{t\in[0,T]}    \|t^{\frac{|s|} {4} } p(\cdot,t)
\| _{L^2(\mathbb{R})}
%
%
\displaystyle \leq   \sup_{t\in[0,T]} t^{\frac{|s|} {4} } \int_0^T
\left\| \chi (\tau)
W_{\mathbb{R}} (t-\tau) f(x, \tau)  \right\| _{L^2 (\mathbb{R})}   d\tau    \\
\ns
\displaystyle \leq  2^{ \frac{|s|}{4} -1 }  \int_0^T
\sup_{t\in[0,T]}   (t-\tau) ^{\frac{|s|} {4} } \left\| \chi (\tau)
W_{\mathbb{R}} (t-\tau) f(x, \tau)  \right\| _{L^2 (\mathbb{R})}   d\tau    \\
\ns
\displaystyle \quad  +   2^{ \frac{|s|}{4} -1 }  \sup_{t\in[0,T]}
\int_0^T \tau^{\frac{|s|} {4} } \left\| \chi (\tau)
W_{\mathbb{R}} (t-\tau) f(x, \tau)  \right\| _{L^2 (\mathbb{R})}   d\tau    \\
\ns
\displaystyle \leq  2^{ \frac{|s|}{4} -1 }    \int_0^T
\sup_{t\in[\tau,T]}  (t-\tau) ^{\frac{|s|} {4} }   \left\|
W_{\mathbb{R}} (t-\tau) f(x, \tau)  \right\| _{L^2 (\mathbb{R})}   d\tau    \\
\ns
\displaystyle \quad  +   2^{ \frac{|s|}{4} -1 }  \sup_{t\in[0,T]}
\int_0^T \tau^{\frac{|s|}{4}}
 \left( \int_{-\infty}^{+\infty}
 \left| e^{ ( i \delta \xi^3 - \xi^4 ) (t-\tau)}   \widehat{f} (\xi,
\tau) \right|^2
  d\xi  \right)^{\frac{1}{2}}   d\tau .
\end{array}
\end{equation}
Then (\ref{2.4.21}) and Proposition \ref{Pro 6.1}  show that
\begin{equation*}
\begin{array}{ll}
\displaystyle   \sup_{t\in[0,T]}   \| t^{\frac{|s|} {4} } p(\cdot,t)
\| _{L^2(\mathbb{R})}\\ \ns
\displaystyle \leq  2^{ \frac{5|s|}{4} -1 }  \int_0^T   \left\|
f(\cdot, \tau)  \right\| _{H^s (\mathbb{R})}   d\tau
 +  2^{ \frac{|s|}{4} -1 } \int_0^T   \tau^{\frac{|s|}{4}}
 \left( \int_{-\infty}^{+\infty}
 \left| \widehat{f} (\xi,
\tau) \right|^2
  d\xi  \right)^{\frac{1}{2}}   d\tau \\ \ns
\leq  2^{ \frac{5|s|}{4} -1 }   \left( \|f \|_{L^1(0,T;
H^s(\mathbb{R}))}   + \|\tau ^{\frac{|s|} {4} } f \|_{L^1 (0,T;
L^2(\mathbb{R}))} \right),
\end{array}
\end{equation*}
which implies the desired estimate (\ref{2.4.20}).

\medskip

{\it Step 3.} We claim that
\begin{equation}\label{2.4.23}
 \| t^{\frac{|s|} {4} }  p\| _{L^2( 0,T; H^{2}
(\mathbb{R}) )}
  \leq   C  2^{ \frac{|s|}{4} }  \left( \|f \|_{L^1(0,T; H^s(\mathbb{R}))}   +  \|t^{\frac{|s|} {4} } f \|_{L^1
(0,T; L^2(\mathbb{R}))} \right).
\end{equation}
Indeed, by (\ref{2.4.15}) and Minkowski's inequality, we see that
\begin{equation}\label{2.4.25}
\begin{array}{ll}
\displaystyle   \|   t^{\frac{|s|} {4} }   p \| _{L^2(0,T; H^2
(\mathbb{R}))}
\displaystyle \leq   \int_0^T \left\|  t^{\frac{|s|} {4} }   \chi
(\tau) W_{\mathbb{R}} (t-\tau) f(x, \tau)  \right\| _{L^2(0,T; H^2
(\mathbb{R}))}   d\tau    \\
\ns
\displaystyle \leq   2^{ \frac{|s|}{4} -1 }  \int_0^T   \left\|
(t-\tau) ^{\frac{|s|} {4} }   \chi (\tau) W_{\mathbb{R}} (t-\tau)
f(x, \tau) \right\| _{L^2(0,T; H^2
(\mathbb{R}))}    d\tau    \\
\ns
\displaystyle \quad  +   2^{ \frac{|s|}{4} -1 }  \int_0^T
 \left\|   \tau^{\frac{|s|} {4} }    \chi (\tau) W_{\mathbb{R}} (t-\tau) f(x,
\tau)  \right\|_{L^2(0,T; H^2
(\mathbb{R}))}   d\tau    \\
\ns
\displaystyle \leq    2^{ \frac{|s|}{4} -1 }    \int_0^T  \left\|
(t-\tau) ^{\frac{|s|} {4} }
 W_{\mathbb{R}} (t-\tau) f(x, \tau)  \right\|
_{L^2(\tau,T; H^2
(\mathbb{R}))}    d\tau   \\
\ns
\displaystyle \quad  +    2^{ \frac{|s|}{4} -1 }  \int_0^T \left\|
\tau^{\frac{|s|}{4}}  W_{\mathbb{R}} (t-\tau) f(x, \tau) \right\|
_{L^2(\tau,T; H^{2} (\mathbb{R}) )} d\tau .
\end{array}
\end{equation}
Note that
\begin{equation}\label{2.4.27}
\begin{array}{ll}
\displaystyle    \left\| \tau^{\frac{|s|}{4}}  W_{\mathbb{R}}
(t-\tau) f(x, \tau) \right\|^2 _{L^2(\tau,T; H^{2} (\mathbb{R}) )}\\
\ns
\displaystyle =   \int_{\tau}^T \int_{-\infty}^{+\infty}
(1+\xi^2)^{2}
 \left| e^{ ( i \delta \xi^3 - \xi^4 ) (t-\tau)}   \tau^{\frac{|s|}{4}}  \widehat{f} (\xi,
\tau) \right|^2
  d\xi  dt  \\ \ns
\displaystyle  =  \int_{-\infty}^{+\infty}  \left(  \int_{\tau}^T
 (1+\xi^2)^{2}  e^{- 2\xi^4(t-\tau) }    dt \right)
\left| \tau^{\frac{|s|}{4}}  \widehat{f} (\xi, \tau) \right|^2  d\xi
\\ \ns
\displaystyle  \leq   \int_{-\infty}^{+\infty}  \left( \int_{\tau}^T
\left( 2+2\xi^4  e^{- 2\xi^4(t-\tau) }  \right)  dt \right)
\left| \tau^{\frac{|s|}{4}}  \widehat{f} (\xi, \tau) \right|^2  d\xi \\
\ns
\leq (2T+1) \left\|\tau^{\frac{|s|}{4}}  f(\cdot,\tau) \right\|^2
_{L^2 (\mathbb{R})}.
\end{array}
\end{equation}
From  (\ref{2.4.25}),  (\ref{2.4.27}) and Proposition \ref{Pro 6.1},
one gets
\begin{equation*}
\begin{array}{ll}
\displaystyle    \|   t^{\frac{|s|} {4} }   p \| _{L^2(0,T; H^2
(\mathbb{R}))}  \\ \ns
\displaystyle \leq   2^{ \frac{ |s|}{4} -1 }    C  \int_0^T \left\|
f(\cdot, \tau)  \right\| _{H^s (\mathbb{R})}   d\tau
 +   2^{ \frac{|s|}{4} -1 }    (2T+1)^{\frac{1}{2}}  \int_0^T    \left\| \tau^{\frac{|s|}{4}}  f(\cdot,\tau) \right\| _{L^2 (\mathbb{R})}  d\tau  \\ \ns
\leq  2^{ \frac{|s|}{4} -1 }  C  \|f \|_{L^1(0,T; H^s(\mathbb{R}))}
+ 2^{ \frac{|s|}{4} -1 }   (2T+1)^{\frac{1}{2}} \|\tau^{\frac{|s|}
{4} } f \|_{L^1 (0,T; L^2(\mathbb{R}))} ,
\end{array}
\end{equation*}
which implies the desired estimate (\ref{2.4.23}).
This completes the proof of Proposition \ref{Pro 6.2}.\endpf


\section{The linear problem on a half line}\label{Sec 3}

In this section, we consider the smoothing properties of the
associated linear equation.

\subsection{Results on the linear problem}\label{Sec 3.1}

To state our  results clearly and concisely, we introduce the
following notations. For any given $a, b\in \mathbb{R}$, set
\begin{equation}\label{2.1*}
\begin{array}{ll}
\vec{h} (t) = (h_1(t), \;  h_2(t)), \qquad
\displaystyle   {\cal H}^{s} (a,b) = H^{\frac{s}{4}+\frac{3}{8}}
(a,b)
\times H^{\frac{s}{4} +\frac{1}{8}} (a,b),\\
\ns
\displaystyle \| \vec{h}  \|_{{\cal H}^{s} (a,b)}  =  \| h_1 \| _{H
^{\frac{s}{4}+\frac{3}{8}}  (a,b) } +  \| h_2 \| _{H
^{\frac{s}{4}+\frac{1}{8}} (a,b) } .
\end{array}
\end{equation}
For any $s\in \mathbb{R}$, $T>0$ and small $\varepsilon > 0$, put
\begin{equation}\label{2.3*}
\begin{array}{ll}
\displaystyle  X_{s,T} =  \Big \{ w\in  C (  [0, T] ; H^s
(\mathbb{R}^+) ) \cap L^2 ( 0,T; H^{s+2} (\mathbb{R}^+)  )  \ \big|
\\  \ns
\displaystyle   \qquad \qquad  \quad  \partial^k_x w \in C (
[0,+\infty); H^{\frac{s-k}{4} +
\frac{3}{8}} (0,T) ), \  k=0,1   \Big \},  \\
\ns
\displaystyle \|w \|_{X_{s,T}} =  \| w \| _{C ( [0,T]; H^s
(\mathbb{R}^+) ) }    +   \| w \|  _{L^2 (0,T; H^{s+2}
(\mathbb{R}^+) )}     \\
\ns
 \displaystyle \qquad \qquad \quad   +  \| w \|  _{ C ( [0,+\infty); H^{\frac{s}{4} + \frac{3}{8}}
(0,T ) )}      +  \| w_x \|  _{ C ( [0,+\infty); H^{\frac{s}{4} +
\frac{1}{8}} (0,T ) )} .
\end{array}
\end{equation}
\begin{equation}\label{6.1.1}
\begin{array}{ll}
\displaystyle  X_{s,T}^{\varepsilon} =  \left \{ w\in C (  [0, T] ;
H^s (\mathbb{R}^+ ) )  \cap   L^2(0,T; H^{s+2}(\mathbb{R}^+ ) )   \
\big| \
 \| w \| _{X_{s,T}^{\varepsilon}} < +\infty \right \},  \\
\ns
\displaystyle \|w \|_{X_{s,T}^{\varepsilon}}  =    \| w \| _{ C (
[0, T] ; H^s (\mathbb{R}^+ ) ) } +  \|w\| _{L^2(0,T;
H^{s+2}(\mathbb{R}^+ ) )}
\\ \ns
\qquad \qquad \quad +  \| t ^{\frac{|s|}{4} + \varepsilon}  w \| _{
C ( [0, T] ; L^2 (\mathbb{R}^+ ) ) }  +
 \| t ^{\frac{|s|}{4} + \varepsilon} w\| _{L^2(0,T; H^{2}(\mathbb{R}^+ ) )} .
\end{array}
\end{equation}

Our results are as follows.

\begin{proposition}\label{Pro 2.1*}
Let $ - 2 \leq s \leq 0$, $T>0$ and  $\vec{h} \in   {\cal H} ^s (
\mathbb{R}^+ )$. Then equation
\begin{eqnarray}\label{2.1.1*}
\left\{\begin{array}{ll} \displaystyle v_{t} + v_{xxxx} + \delta
v_{xxx} = 0, & (x,t)\in \mathbb{R}^+ \times \mathbb{R}^+,
\\ \ns
v(x,0)= 0, &  x\in \mathbb{R}^+,
\\ \ns
v(0,t)=h_1(t), \quad v_x(0,t)=h_2(t),  \qquad  & t\in \mathbb{R}^+,
\end{array}\right.
\end{eqnarray}
admits a unique solution $v  \in C ( [0,+\infty) ; H^s
(\mathbb{R}^+) ) \cap L^2  (  \mathbb{R}^+  ; H^{s+2} (\mathbb{R}^+)
)$ with
$\partial^k_x v \in$ \linebreak  $ C ( [0,+\infty); H^{\frac{s-k}{4}
+ \frac{3}{8}} (\mathbb{R}^+ ) )$ ($k=0,1$).
Moreover, there exists a constant $C> 0$ such that for any $\vec{h}
\in {\cal H} ^s (\mathbb{R}^+)$, it holds
\begin{equation}\label{2.1.2*}
\|v \|_{X_{s,T}}   \leq C \| \vec{h} \| _{{\cal H} ^s
(\mathbb{R}^+)}.
\end{equation}
Furthermore, if $-2  < s < 0$,  $\varepsilon>0$  and
$t^{\frac{|s|}{4} + \varepsilon} \vec{h} \in {\cal H} ^{0} (0,T)$,
it holds
\begin{equation}\label{2.1.3*}
\| v \|_{X_{s,T}^{\varepsilon} }  \leq C  \left ( \| \vec{h} \|_{
{\cal H}^s (0,T) } + \| t^{\frac{|s|} {4} + \varepsilon } \vec{h} \|
_{ {\cal H}^0 (0,T) } \right).
\end{equation}

\end{proposition}

\begin{proposition}\label{Pro 2.2*}
Let $- 2\leq s \leq 0$, $T>0$ and $\phi \in  H^s (\mathbb{R}^+) $.
Then equation
\begin{eqnarray}\label{2.2.1*}
\left\{\begin{array}{ll} \displaystyle v_{t} + v_{xxxx} + \delta
v_{xxx} = 0, & (x,t)\in \mathbb{R}^+ \times (0,T),
\\ \ns
v(x,0)= \phi (x), &  x\in \mathbb{R}^+,
\\ \ns
v(0,t)= 0, \quad v_x(0,t)= 0,  \qquad  & t\in (0,T),
\end{array}\right.
\end{eqnarray}
admits a unique solution $v  \in C ( [0,T] ; H^s (\mathbb{R}^+) )
\cap L^2  ( 0,T ; H^{s+2} (\mathbb{R}^+)  )$ with
$\partial^k_x v \in $  \linebreak   $C (   [0, +\infty);
H^{\frac{s-k}{4} + \frac{3}{8}} (0,T) )$ ($k=0,1$).
Moreover, there exist a constant $C> 0$ such that for any $\phi \in
H^s(\mathbb{R})$, it holds
\begin{equation}\label{2.2.2*}
\|v \|_{X_{s,T}}  \leq   C  \| \phi \| _{H^s (\mathbb{R}^+)}.
\end{equation}
Furthermore, if $-2  < s < 0$, $0<T\leq 1$ and  $\varepsilon>0$,  it
holds
\begin{equation}\label{2.2.3*}
\| v \|_{X_{s,T}^{\varepsilon} }  \leq C \| \phi \| _{H
^s(\mathbb{R}^+)}.
\end{equation}

\end{proposition}

\begin{proposition}\label{Pro 2.3*}
Let $-2 \leq s \leq 0$, $T>0$ and $f \in L^1(0,T;
H^s(\mathbb{R}^+))$. Then equation
\begin{eqnarray}\label{2.3.1*}
\left\{\begin{array}{ll} \displaystyle v_{t} + v_{xxxx} + \delta
v_{xxx} = f(x,t), & (x,t)\in \mathbb{R}^+ \times (0,T),
\\ \ns
v(x,0)= 0, &  x\in \mathbb{R}^+,
\\ \ns
v(0,t)= 0, \quad v_x(0,t)= 0,  \qquad  & t\in (0,T),
\end{array}\right.
\end{eqnarray}
admits a unique solution $v  \in C ( [0,T] ; H^s (\mathbb{R}^+) )
\cap L^2  ( 0,T ; H^{s+2} (\mathbb{R}^+)  )$ with
$\partial^k_x v \in $  \linebreak   $C (   [0, +\infty);
H^{\frac{s-k}{4} + \frac{3}{8}} (0,T) )$ ($k=0,1$).
Moreover, there exist a constant $C> 0$ such that for any $f \in
L^1(0,T; H^s(\mathbb{R}^+))$, it holds
\begin{equation}\label{2.3.2*}
\|v \|_{X_{s,T}}     \leq C   \| f \| _{L^1(0,T; H
^s(\mathbb{R}^+))}.
\end{equation}
Furthermore, if $-2  < s < 0$,  $0<T\leq 1$,  $\varepsilon>0$  and $
t^{\frac{|s|} {4}  +  \varepsilon }  f  \in L^1 (0,T; L^2 (
\mathbb{R}^+ ) )$, it holds
\begin{equation}\label{2.3.3*}
\| v \|_{X_{s,T}^{\varepsilon} } \leq C \left( \|f \|_{L^1(0,T;
H^s(\mathbb{R}^+ ))} + \|t^{\frac{|s|} {4}  +   \varepsilon } f
\|_{L^1 (0,T; L^2(\mathbb{R}^+ ))}   \right).
\end{equation}

\end{proposition}

We will prove estimates (\ref{2.1.2*}), (\ref{2.2.2*}) and
(\ref{2.3.2*}) in Subsection \ref{Sec 3.2}. Estimates (\ref{2.1.3*})
will be proven in Subsection \ref{Sec 3.3}.
We omit the proof of  (\ref{2.2.3*}) and (\ref{2.3.3*}) since they
are similar to the proof of  (\ref{2.2.2*}) and Proposition \ref{Pro
6.2}, respectively.


\subsection{Proof of estimates (\ref{2.1.2*}), (\ref{2.2.2*}) and (\ref{2.3.2*})}\label{Sec 3.2}

We recall the following lemma, which follows from minor
modifications of  Lemma 3.1 in \cite{Zhang1}.

\begin{lemma}\label{Lma 2.1.1*}
Let $\gamma (\rho)$ be a continuous complex-valued function defined
on $(0,+\infty)$ satisfying the following three conditions:
\begin{enumerate}
\item[i)] $\displaystyle {\text Re} \gamma(\rho) <0$, for $\rho >0$;

\item[ii)]There exist $\delta>0$ and $b>0$ such that $\displaystyle \sup_{0<\rho <\delta }  \frac{|{\text Re} \gamma(\rho)|}  {\rho}   \geq b$;

\item[iii)]There exists a complex number $\alpha+i \beta $ with $\alpha<0$ such that $\displaystyle \lim _{\rho\rightarrow +\infty} \frac{|\gamma(\rho)|}  {\rho}   =\alpha +i \beta$.
\end{enumerate}
Then there exists a constant $C>0$ such that for all $f\in L^2(0,
+\infty)$,
$$ \left \| \int_0^{+\infty} e ^{\gamma (\rho) x} f(\rho) d \rho \right \| _{L^2(\mathbb{R}^+)}  \leq C  \|f(\cdot)\|_{L^2(\mathbb{R}^+)}.$$
\end{lemma}

\noindent {\bf Proof of estimate (\ref{2.1.2*}). }
We divide the proof into four steps. In Step 1, an explicit solution
formula of equation (\ref{2.1.1*}) is given. Then with the help of
this solution formula, estimate (\ref{2.1.2*}) is proved for $-2
\leq s\leq 0$ in Step 2 - Step 4.

\medskip

{\it Step 1.}
Applying the Laplace transform with respect to $t$ in equation
(\ref{2.1.1*}), we see that for any $\tau$ with $\text{Re}\tau>0$,
\begin{eqnarray}\label{2.1.4*}
\left\{\begin{array}{ll} \displaystyle \tau\widetilde{v} +
\widetilde{v}_{xxxx}  + \delta \widetilde{v}_{xxx} =0,  &
(x,\tau)\in \mathbb{R}^+ \times \mathbb{R}^+,
\\ \ns
\widetilde{v}(0,\tau)=\widetilde{h}_1(\tau), \quad
\widetilde{v}_x(0,\tau)=\widetilde{h}_2(\tau),  \qquad &  \tau \in
\mathbb{R}^+,
\end{array}\right.
\end{eqnarray}
where
$$\widetilde{v}(x,\tau)=\int_0^{+\infty} e^{-\tau t} v(x,t) dt,\qquad\quad
\widetilde{h}_j(\tau)=\int_0^{+\infty} e^{-\tau t} h_j(t) dt, \quad
j=1,2.$$
The solution of equation (\ref{2.1.4*}) can be expressed as follows
\begin{equation}\label{2.1.5*}
\widetilde{v} (x,\tau) = c_1(\tau) e^{\lambda _1 (\tau)x}  +
c_2(\tau) e^{\lambda _2 (\tau)x}.
\end{equation}
Here, $\lambda_1(\tau)$ and $\lambda_2(\tau)$ are the solutions of
the characteristic equation
\begin{equation}\label{2.1.6*}
\lambda ^4 + \delta \lambda^3  + \tau = 0
\end{equation}
with negative real parts,
and $c_1 =c_1 (\tau)$, $c_2 =c_2 (\tau)$ are the solutions of
\begin{eqnarray}\label{2.1.7*}
\left\{\begin{array}{ll}
\displaystyle c_1  + c_2  = \widetilde{h}_1 (\tau),\\
\ns\ns \displaystyle c_1 \lambda_1 (\tau)   + c_2 \lambda_2 (\tau)
  = \widetilde{h}_2 (\tau).
\end{array}\right.
\end{eqnarray}
By (\ref{2.1.5*})-(\ref{2.1.7*}),  for any fixed $r>0$, the solution
of (\ref{2.1.1*}) can be represented in the form
\begin{equation}\label{2.1.8*}
\begin{array}{ll}
v(x,t)  &  \displaystyle = \frac{1}{2\pi i}\int_{r-i \infty}^{r+i
\infty} e^{\tau t} \widetilde{v}(x,\tau) d\tau \\ \ns
& \displaystyle = \frac{1}{2\pi i}\int_{r-i \infty}^{r+i \infty}
e^{\tau t} \left( \frac{ \lambda_2 \widetilde{h}_1 (\tau) -
\widetilde{h}_2 (\tau) }{\lambda_2 - \lambda_1} e^{\lambda _1  x}
  +
\frac{ \widetilde{h}_2 (\tau)  - \lambda_1 \widetilde{h}_1 (\tau)
}{\lambda_2 - \lambda_1} e^{\lambda _2  x}   \right)
   d\tau.
\end{array}
\end{equation}
Since the right-hand side of (\ref{2.1.8*}) is continuous with
respect to $r$ and the left-hand side does not depend on $r$, we can
take $r=0$ in (\ref{2.1.8*}). Setting $\tau= i 8 \rho^4 $ with
$0\leq \rho <+\infty$ in (\ref{2.1.6*}), then the two solutions of
characteristic equation (\ref{2.1.6*}) with $\delta =0 $ and with
negative real parts are as follows
\begin{equation*}
- \rho \sqrt{\sqrt{2}+1}
  +i \rho \sqrt{\sqrt{2}-1},  \quad
 - \rho \sqrt{\sqrt{2}-1}
  -i \rho \sqrt{\sqrt{2}+1}.
\end{equation*}
Thus, as $\rho \rightarrow + \infty$, we have the asymptotic
expression of the two solutions of characteristic equation
(\ref{2.1.6*}) with negative real parts:
\begin{equation}\label{2.1.10*}
\lambda_1^+ (\rho) \sim - \rho \sqrt{\sqrt{2}+1}
  +i \rho \sqrt{\sqrt{2}-1},  \quad
\lambda_2^+ (\rho) \sim - \rho \sqrt{\sqrt{2}-1}
  -i \rho \sqrt{\sqrt{2}+1}.
\end{equation}
By (\ref{2.1.8*}) and (\ref{2.1.10*}), it is easy to check that
\begin{equation}\label{2.1.15*}
v(x,t) = v^+(x,t)  +  \overline{ v^+(x,t)},
\end{equation}
where $ \overline{ v^+(x,t)}$ is the conjugation of  $ v^+(x,t)$ and
\begin{equation}\label{2.1.12*}
\begin{array}{rr}
 v^+(x,t)  =
\displaystyle \frac{16}{\pi}\int_0^{+ \infty} e^{ i 8\rho^4 t }
\left( e^{\lambda _1^+  x} \frac{ \lambda_2 ^+ }{\lambda_2 ^+ -
\lambda_1^+ }  \widetilde{h}_1 ( i 8 \rho^4 )
- e^{\lambda _2^+  x} \frac{ \lambda_1 ^+}{\lambda_2 ^+ -
\lambda_1^+ } \widetilde{h}_1 ( i 8 \rho^4 )  \right)\rho^3 d\rho\\
\ns
 \displaystyle + \frac{16}{\pi}\int_0^{+ \infty} e^{ i 8\rho^4 t }
\left( - e^{\lambda _1^+  x} \frac{ 1 }{\lambda_2 ^+  - \lambda_1^+
} \widetilde{h}_2 ( i 8 \rho^4 )
+ e^{\lambda _2^+  x} \frac{ 1 }{\lambda_2 ^+ - \lambda_1^+}
\widetilde{h}_2 ( i 8 \rho^4 )  \right)\rho^3 d\rho.
\end{array}
\end{equation}

{\it Step 2.} We claim that for $-2 \leq s\leq 0$, it holds
\begin{equation}\label{2.1.30*}
\sup _{ t\in [0,T] } \|v (\cdot,t)\| _{H^s(\mathbb{R}^+)}
\leq C \left( \| h_1 \| _{ H^{\frac{s}{4}+\frac{3}{8}}
(\mathbb{R}^+)}  + \| h_2 \| _{ H^{\frac{s}{4}+\frac{1}{8}}
(\mathbb{R}^+)}  \right).
\end{equation}
It suffices to prove that estimate (\ref{2.1.30*}) is true for
$s=0,-1,-2$.  For non-integer negative values of $s$, estimate
(\ref{2.1.30*}) can be obtained by standard interpolation theory.

Now  we consider the case $s=0,-1,-2$.
For any $ \varphi (x) \in H_0^{-s} (\mathbb{R}^+)$, we integrate by
parts and deduce that
\begin{equation*}
\begin{array}{ll}
\displaystyle  \left \langle      \int _0 ^{+\infty}  e^{\lambda^+
_j (\rho) x } F(\rho,t) d\rho,     \varphi (x) \right\rangle
_{H^s(\mathbb{R}^+) \times H_0^{-s}(\mathbb{R}^+) }\\ \ns
  \displaystyle   = \int _0^{+ \infty}  \left( \int _0 ^{+\infty} e^{\lambda
_j^+ (\rho) x } F(\rho, t) d\rho \right)  \varphi (x)  dx\\ \ns
 \displaystyle  =  \int _0^{+ \infty}  (-1)^s  \left(
   \int _0^{+\infty}  \left( \lambda _j^+  (\rho) \right) ^s   e^{\lambda _j^+  (\rho) x }  F(\rho, t) d\rho \right)  \varphi^{(-s)}(x)
   dx,
\end{array}
\end{equation*}
which implies
\begin{equation}\label{2.1.31}
\begin{array}{ll}
\displaystyle    \left \| \int _0 ^{+\infty}  e^{\lambda^+ _j (\rho)
x } F(\rho,t) d\rho  \right \| ^2_{H^s(\mathbb{R}^+) }
\leq \left\|
   \int _0^{+\infty}  \left( \lambda _j^+  (\rho) \right) ^s    e^{\lambda _j^+  (\rho) x }  F(\rho,t) d\rho
   \right\|^2  _{L^2(\mathbb{R}^+)}.
\end{array}
\end{equation}
From (\ref{2.1.12*}) and (\ref{2.1.31}), we see that for any $t\in
 [0,T] $,
\begin{equation}\label{2.1.32*}
\begin{array}{ll}
\displaystyle  \|  v^+ (\cdot,t) \| _{H^s(\mathbb{R}^+)} ^2
 &  \displaystyle \leq C   \left\|  \int_0^{+ \infty}    \left( \lambda _1^+   \right) ^s  e^{ i 8\rho^4 t }
 e^{\lambda _1^+  x} \frac{ \lambda_2 ^+ }{\lambda_2^+  - \lambda_1^+ }
\widetilde{h}_1 ( i 8 \rho^4 ) \rho^3 d\rho   \right\|
_{L^2(\mathbb{R}^+)} ^2 \\ \ns
& \displaystyle  \quad + C   \left\|  \int_0^{+ \infty} -   \left(
\lambda _2^+   \right) ^s   e^{ i 8\rho^4 t }
 e^{\lambda _2^+  x} \frac{ \lambda_1^+ }{\lambda_2^+  - \lambda_1^+ }
\widetilde{h}_1 ( i 8 \rho^4 )  \rho^3 d\rho   \right\|
_{L^2(\mathbb{R}^+)} ^2 \\ \ns
& \displaystyle \quad + C   \left\|  \int_0^{+ \infty} -  \left(
\lambda _1^+   \right) ^s    e^{ i 8\rho^4 t }
 e^{\lambda _1^+  x} \frac{ 1}{\lambda_2^+ - \lambda_1^+}
\widetilde{h}_2 ( i 8 \rho^4 )   \rho^3 d\rho   \right\|
_{L^2(\mathbb{R}^+)} ^2 \\ \ns
& \displaystyle  \quad + C   \left\|  \int_0^{+ \infty}    \left(
\lambda _2^+   \right) ^s     e^{ i 8\rho^4 t }
 e^{\lambda _2^+  x} \frac{ 1 }{\lambda_2^+ - \lambda_1^+}
\widetilde{h}_2 ( i 8 \rho^4 )
\rho^3 d\rho   \right\| _{L^2(\mathbb{R}^+)} ^2.
\end{array}
\end{equation}
It follows from (\ref{2.1.32*}) and Lemma \ref{Lma 2.1.1*} that for
any $t\in  [0,T]  $,
\begin{equation*}
\begin{array}{ll}
\displaystyle  \|  v^+ (\cdot,t) \| _{H^s (\mathbb{R}^+)} ^2 \\ \ns
\displaystyle \leq C    \int_0^{+ \infty}    \left|  \left( \lambda
_1^+   \right) ^{s}   \frac{ \lambda_2^+ }{\lambda_2^+  -
\lambda_1^+ } \widetilde{h}_1 ( i 8 \rho^4 )  \rho^3   \right|^2
d\rho
\displaystyle + C   \int_0^{+ \infty} \left|  \left( \lambda _2^+
  \right) ^{s}   \frac{ \lambda_1^+ }{\lambda_2^+  -
\lambda_1^+ }   \widetilde{h}_1 ( i 8 \rho^4 )  \rho^3  \right|^2   d\rho \\
\ns
\displaystyle  \quad  + C    \int_0^{+ \infty} \left|  \left(
\lambda _1^+  \right) ^s    \frac{ 1}{\lambda_2^+  - \lambda_1^+ }
\widetilde{h}_2 ( i 8 \rho^4 )   \rho^3  \right|^2  d\rho
\displaystyle    + C   \int_0^{+ \infty}  \left|  \left( \lambda
_2^+  \right) ^s    \frac{ 1 }{\lambda_2^+  - \lambda_1^+ }
\widetilde{h}_2 ( i 8 \rho^4 )   \rho^3 \right|^2  d\rho.
\end{array}
\end{equation*}
Combining  (\ref{2.1.10*}), the above inequality gives
\begin{equation}\label{2.1.34*}
\begin{array}{ll} \displaystyle \sup_{t\in [0,T] }   \|  v^+ (\cdot,t) \| _{H^s (\mathbb{R}^+)} ^2
 \displaystyle  \leq C
 \int_0^{+ \infty}    \rho ^{2s + 6}
\left| \widetilde{h}_1 (i 8 \rho^4 ) \right| ^2  d\rho  +  C
 \int_0^{+ \infty}   \rho ^{2s + 4}
\left| \widetilde{h}_2 (i 8 \rho^4 ) \right| ^2  d\rho.
\end{array}
\end{equation}
Setting $\mu = 8 \rho^4$ in (\ref{2.1.34*}) and recalling
$\widetilde{h}_j(\tau)=\int_0^{+\infty} e^{-\tau t} h_j(t) dt$, we
deduce that
\begin{equation*}\label{2.1.36}
\begin{array}{ll}
\displaystyle \sup_{t\in [0,T] } \| v^+ (\cdot,t) \|
_{H^s(\mathbb{R}^+)} ^2
& \displaystyle    \leq C \int_0^{+ \infty} (1+
\mu^2)^{\frac{s}{4}+\frac{3}{8}}
 \left| \int_0^{+ \infty} e^{-i\mu \tau}  h_1 (\tau) d\tau \right| ^2
 d\mu   \\ \ns
& \displaystyle \quad + C \int_0^{+ \infty} (1+
\mu^2)^{\frac{s}{4}+\frac{1}{8}}
 \left| \int_0^{+ \infty} e^{-i\mu \tau}  h_2 (\tau) d\tau \right| ^2  d\mu\\ \ns
& \displaystyle \leq C  \left( \| h_1 \|^2 _{
H^{\frac{s}{4}+\frac{3}{8}} (\mathbb{R}^+)}  + \| h_2 \|^2 _{
H^{\frac{s}{4}+\frac{1}{8}} (\mathbb{R}^+)}  \right),
\end{array}
\end{equation*}
which combining (\ref{2.1.15*}) yields that (\ref{2.1.30*}) holds
for $s=0,-1,-2$.

\medskip

{\it Step 3.} We claim that for $-2 \leq s\leq 0$, it holds
\begin{equation}\label{2.1.38*}
\|v\|_{ L^2 ( \mathbb{R}^+ ; H^{s+2} (\mathbb{R}^+)  )}
\leq C \left( \| h_1 \| _{ H^{\frac{s}{4}+\frac{3}{8}}
(\mathbb{R}^+)}  + \| h_2 \| _{ H^{\frac{s}{4}+\frac{1}{8}}
(\mathbb{R}^+)}  \right).
\end{equation}
By standard interpolation theory, it suffices to show that
(\ref{2.1.38*}) is true  for $s=0,-1,-2$.

Setting $\mu = 8 \rho^4$ in (\ref{2.1.12*}), we have that for
$k=0,1,2$,
\begin{equation}\label{2.1.44*}
\begin{array}{ll}
\partial^{k} _x v^+ (x,t)
& \displaystyle  = \frac{1}{2\pi}\int_0^{+ \infty} e^{  i \mu t }
e^{\lambda _1^+  x} \frac{ \lambda_2^+ }{\lambda_2^+  - \lambda_1^+
} \left( \lambda _1^+  \right)^{k}
 \widetilde{h}_1 (i \mu)d\mu  \\ \ns
&  \displaystyle \quad   -  \frac{1}{2\pi}\int_0^{+ \infty} e^{  i
\mu t }    e^{\lambda _2^+  x} \frac{ \lambda_1^+ }{\lambda_2^+  -
\lambda_1^+ }       \left( \lambda _2^+  \right)^{k}
 \widetilde{h}_1 (i \mu)d\mu  \\ \ns
&  \displaystyle \quad   -  \frac{1}{2\pi}\int_0^{+ \infty} e^{  i
\mu t }    e^{\lambda _1^+   x} \frac{ 1 }{\lambda_2^+  -
\lambda_1^+ } \left( \lambda _1^+  \right)^{k}
 \widetilde{h}_2 (i \mu)d\mu  \\ \ns
&  \displaystyle \quad   +  \frac{1}{2\pi}\int_0^{+ \infty} e^{  i
\mu t }    e^{\lambda _2^+  x} \frac{ 1 }{\lambda_2^+  - \lambda_1^+
} \left( \lambda _2^+  \right)^{k}
 \widetilde{h}_2 (i \mu)d\mu.
\end{array}
\end{equation}
By (\ref{2.1.44*}) and Plancherel's theorem, one gets that for any
$x \in \mathbb{R}^+$,
\begin{equation*}
\begin{array}{ll}
\displaystyle  \| \partial ^{k}_x v^+ (x, \cdot) \| _{L^2
(\mathbb{R}^+)} ^2
& \displaystyle  \leq  C   \int_0^{+ \infty}
 \left| e^{\lambda _1^+  x}
  \frac{ \lambda_2^+ }{\lambda_2^+  - \lambda_1^+ }
  \left(\lambda _1^+ \right)^{k} \widetilde{h}_1 (i\mu) \right| ^2 d\mu  \\ \ns
& \displaystyle  \quad  +  C  \int_0^{+ \infty}
 \left| e^{\lambda _2^+  x}
  \frac{ \lambda_1^+ }{\lambda_2^+  - \lambda_1^+ }
  \left(\lambda _2^+  \right)^{k} \widetilde{h}_1
(i\mu) \right| ^2 d\mu  \\ \ns
& \displaystyle  \quad   +  C \int_0^{+ \infty}
 \left| e^{\lambda _1^+  x}
  \frac{ 1 }{\lambda_2^+  - \lambda_1^+ }
  \left(\lambda _1^+  \right)^{k} \widetilde{h}_2
(i\mu) \right| ^2 d\mu  \\ \ns
&  \displaystyle  \quad +  C  \int_0^{+ \infty}
 \left| e^{\lambda _2^+  x}
  \frac{ 1 }{\lambda_2^+  - \lambda_1^+ }
  \left(\lambda _2^+ \right)^{k} \widetilde{h}_2
(i\mu) \right| ^2 d\mu,
\end{array}
\end{equation*}
which gives  that for $k=0,1,2$,
\begin{equation}\label{2.1.39*}
\begin{array}{ll}
\displaystyle  \int_0^{+\infty} \| \partial ^{k}_x v^+ (x,\cdot)
\|^2 _{L^2 ( \mathbb{R}^+ )}  dx  \\
\ns
\displaystyle  \leq  C   \int_0^{+ \infty}
 \left(  \int_0^{+ \infty}  e^{2 \lambda _1^+  x}  dx  \right)
  \left|  \frac{ \lambda_2^+ }{\lambda_2^+  - \lambda_1^+ }
  \left(\lambda _1^+  \right)^{k} \widetilde{h}_1
(i\mu) \right| ^2 d\mu  \\ \ns
\displaystyle  \quad  +  C  \int_0^{+ \infty}
 \left(  \int_0^{+ \infty}  e^{2 \lambda _2^+  x}  dx  \right)
  \left| \frac{ \lambda_1^+ }{\lambda_2^+  - \lambda_1^+ }
  \left(\lambda _2^+  \right)^{k} \widetilde{h}_1
(i\mu) \right| ^2 d\mu  \\ \ns
\displaystyle  \quad   +  C \int_0^{+ \infty}
  \left(  \int_0^{+ \infty}  e^{2 \lambda _1^+  x}  dx  \right)
  \left| \frac{ 1 }{\lambda_2^+  - \lambda_1^+ }
  \left(\lambda _1^+  \right)^{k} \widetilde{h}_2
(i\mu) \right| ^2 d\mu  \\ \ns
\displaystyle  \quad +  C  \int_0^{+ \infty}
 \left(  \int_0^{+ \infty}  e^{2 \lambda _2^+  x}  dx  \right)
  \left|  \frac{ 1 }{\lambda_2^+  - \lambda_1^+ }
  \left(\lambda _2^+  \right)^{k} \widetilde{h}_2
(i\mu) \right| ^2 d\mu.
\end{array}
\end{equation}
Therefore, combining (\ref{2.1.10*}) and (\ref{2.1.39*}), recalling
$\mu = 8 \rho^4$,
\begin{equation}\label{2.1.40*}
\begin{array}{ll}
\displaystyle   \|  v^+  \| _{L^2 ( \mathbb{R}^+ ; H^{s+2}
(\mathbb{R}^+)  )} ^2        = \sum_{k=0}^{s+2}  \| \partial ^{k}_x
v^+ \| _{L^2 ( \mathbb{R}^+ ; L^{2} (\mathbb{R}^+)  )} ^2
 = \sum_{k=0}^{s+2}  \int_0^{+\infty} \| \partial ^{k}_x v^+ (x,\cdot)
\|^2 _{L^2 ( \mathbb{R}^+ )}  dx  \\
\ns
 \displaystyle  \leq C
\sum_{k=0}^{s+2}      \int_0^{+ \infty}   \mu^{2\cdot
{\frac{2k-1}{8}}} \left| \widetilde{h}_1 (i\mu) \right| ^2  d\mu
+ C  \sum_{k=0}^{s+2}   \int_0^{+ \infty}    \mu^{2\cdot
{\frac{2k-3}{8}}}
 \left| \widetilde{h}_2
(i \mu) \right| ^2  d\mu\\ \ns
\displaystyle \leq C  \left( \| h_1 \| ^2 _{
H^{\frac{s}{4}+\frac{3}{8}} (\mathbb{R}^+)}  + \| h_2 \| ^2 _{
H^{\frac{s}{4}+\frac{1}{8}} (\mathbb{R}^+)}  \right).
\end{array}
\end{equation}
From (\ref{2.1.15*}) and (\ref{2.1.40*}), we obtain the desired
estimate (\ref{2.1.38*}) for $s=0,-1,-2$.

\medskip

{\it Step 4.}  We claim that for $-2 \leq s\leq 0$, $\partial^k_x v
\in C ( [0,+ \infty) ; H^{\frac{s-k}{4} + \frac{3}{8}} (\mathbb{R}^+
) )$ $(k=0,1)$. Moreover, it holds
\begin{equation}\label{2.1.46*}
\sup _{x\in  [0, +\infty) } \| \partial ^k _x  v(x, \cdot) \|
_{H^{\frac{s-k}{4} + \frac{3}{8}} (\mathbb{R}^+)} \leq C \left( \|
h_1 \| _{ H^{\frac{s}{4}+\frac{3}{8}} (\mathbb{R}^+)}  + \| h_2 \|
_{ H^{\frac{s}{4}+\frac{1}{8}} (\mathbb{R}^+)}  \right).
\end{equation}

It follows from (\ref{2.1.44*}) that for any $x$, $x_0 \in
[0,+\infty) $,
\begin{equation*}
\begin{array}{ll}
\displaystyle  \left\| \partial^k_x v^+ (x,\cdot)  -  \partial^k_x v^+ (x_0,\cdot) \right\|  _{H^ {\frac{s-k}{4} + \frac{3}{8}} (\mathbb{R}^+) } ^2 \\
\ns
\displaystyle \leq  C
 \int_0^{+ \infty}
(1+ \mu^2)^{\frac{s-k}{4} + \frac{3}{8}} \left|\left( e^{\lambda
_1^+  x} -e^{\lambda _1^+  x_0} \right) \frac{ \lambda_2^+
}{\lambda_2^+  - \lambda_1^+ }
  \left( \lambda _1^+  \right)^k \widetilde{h}_1
(i\mu) \right| ^2 d\mu  \\  \ns
\displaystyle  \quad   + C
 \int_0^{+ \infty}
(1+ \mu^2)^{\frac{s-k}{4} + \frac{3}{8}} \left|\left( e^{\lambda
_2^+  x} -e^{\lambda _2^+  x_0} \right) \frac{ \lambda_1^+
}{\lambda_2^+  - \lambda_1^+ }
  \left( \lambda _2^+  \right)^k \widetilde{h}_1
(i\mu) \right| ^2 d\mu  \\  \ns
\displaystyle \quad  +  C
 \int_0^{+ \infty}
(1+ \mu^2)^{\frac{s-k}{4} + \frac{3}{8}} \left|\left( e^{\lambda
_1^+  x} -e^{\lambda _1^+  x_0} \right) \frac{ 1 }{\lambda_2^+  -
\lambda_1^+ }
  \left( \lambda _1^+  \right)^k \widetilde{h}_2
(i\mu) \right| ^2 d\mu  \\  \ns
\displaystyle   \quad   +   C
 \int_0^{+ \infty}
(1+ \mu^2)^{\frac{s-k}{4} + \frac{3}{8}} \left|\left( e^{\lambda
_2^+  x} -e^{\lambda _2^+  x_0} \right) \frac{ 1 }{\lambda_2^+  -
\lambda_1^+ }
  \left( \lambda _2^+  \right)^k \widetilde{h}_2
(i\mu) \right| ^2 d\mu.
\end{array}
\end{equation*}
From Fatou Lemma, the above inequality shows that
\begin{equation*}
\begin{array}{ll}
\displaystyle  \lim_{x\rightarrow x_0} \left\| \partial^k_x v^+
(x,\cdot)  -  \partial^k_x v^+ (x_0,\cdot) \right\| ^2  _{H^
{\frac{s-k}{4} + \frac{3}{8}} (\mathbb{R}^+) }
=0.
\end{array}
\end{equation*}
Thus,  $\partial^k_x v \in C ( [0,+ \infty) ; H^{\frac{s-k}{4} +
\frac{3}{8}} (\mathbb{R}^+ ) )$ $(k=0,1)$.

It follows from  (\ref{2.1.44*}) that for any $x\in [0, +\infty)$,
\begin{equation*}
\begin{array}{ll}
\displaystyle  \| \partial ^{k}_x v^+ (x, \cdot) \|
_{H^{\frac{s-k}{4} + \frac{3}{8} }(\mathbb{R}^+)} ^2  \\ \ns
\displaystyle  \leq  C   \int_0^{+ \infty}  (1+
\mu^2)^{\frac{s-k}{4} + \frac{3}{8}}
 \left| e^{\lambda _1^+  x}
  \frac{ \lambda_2^+ }{\lambda_2^+ - \lambda_1^+ }
  \left(\lambda _1^+  \right)^{k} \widetilde{h}_1 (i\mu) \right| ^2 d\mu  \\ \ns
\displaystyle  \quad  +  C  \int_0^{+ \infty}  (1+
\mu^2)^{\frac{s-k}{4} + \frac{3}{8}}
 \left| e^{\lambda _2^+  x}
  \frac{ \lambda_1^+ }{\lambda_2^+ - \lambda_1^+ }
  \left(\lambda _2^+  \right)^{k} \widetilde{h}_1
(i\mu) \right| ^2 d\mu  \\ \ns
\displaystyle  \quad   +  C \int_0^{+ \infty}  (1+
\mu^2)^{\frac{s-k}{4} + \frac{3}{8}}
 \left| e^{\lambda _1^+  x}
  \frac{ 1 }{\lambda_2^+ - \lambda_1^+ }
  \left(\lambda _1^+  \right)^{k} \widetilde{h}_2
(i\mu) \right| ^2 d\mu  \\ \ns
\displaystyle  \quad +  C  \int_0^{+ \infty}  (1+
\mu^2)^{\frac{s-k}{4} + \frac{3}{8}}
 \left| e^{\lambda _2^+  x}
  \frac{ 1 }{\lambda_2^+ - \lambda_1^+ }
  \left(\lambda _2^+  \right)^{k} \widetilde{h}_2
(i\mu) \right| ^2 d\mu.
\end{array}
\end{equation*}
Therefore, combining (\ref{2.1.10*}), the above equality shows that
\begin{equation}\label{2.1.48}
\begin{array}{ll} \displaystyle \sup_{x\in [0, +\infty) }\| \partial ^{k}_x v^+ (x, \cdot) \| _{H^ {\frac{s-k}{4} + \frac{3}{8}} (\mathbb{R}^+) } ^2 \\ \ns
 \displaystyle  \leq C
 \int_0^{+ \infty}   (1¡¡+ \mu^2) ^{\frac{s}{4} + \frac{3}{8}}
\left| \widetilde{h}_1 (i\mu) \right| ^2  d\mu
+ C \int_0^{+ \infty}   (1¡¡+ \mu^2) ^{\frac{s}{4} + \frac{1}{8}}
 \left| \widetilde{h}_2
(i \mu) \right| ^2  d\mu\\ \ns
\displaystyle \leq C \left( \| h_1 \| ^2 _{
H^{\frac{s}{4}+\frac{3}{8}} (\mathbb{R}^+)}  + \| h_2 \| ^2 _{
H^{\frac{s}{4}+\frac{1}{8}} (\mathbb{R}^+)}  \right).
\end{array}
\end{equation}
From (\ref{2.1.15*}) and (\ref{2.1.48}), we obtain the desired
estimate (\ref{2.1.46*}) for $-2 \leq s\leq 0$.
\endpf

\begin{remark}\label{Rem 2.1}

Let $W_{bdr} (t):  {\cal H}^{s} (\mathbb{R}^+)   \rightarrow   C (
\mathbb{R}^+_0 ; H^{s} (\mathbb{R}^+) )  \cap L^2 ( \mathbb{R}^+ ;
H^{s+2} (\mathbb{R}^+) ) $ be the solution map of equation
(\ref{2.1.1*}).
Then $W_{bdr} (t) \vec{h}$ can be represented explicitly by
(\ref{2.1.10*})-(\ref{2.1.12*}).

\end{remark}

\noindent {\bf Proof of estimate (\ref{2.2.2*}). }
Let $\bar{\phi} (x)$ be an continuous extension of $\phi (x)$ to the
whole line $\mathbb{R}$ such that
$$\| \bar{\phi} \|_{H^s (\mathbb{R})}  \leq  C \| \phi \| _{H^s(\mathbb{R}^+)} .$$
%
%
%
From $\bar{\phi} \in H^s (\mathbb{R}) $ and Proposition \ref{Pro
6.1}, we have
$$\left(  (W_{\mathbb{R}}(t) \bar{\phi})  (0, t),   \;
(W_{\mathbb{R}}(t) \bar{\phi})_x (0, t) \right)  \in {\cal H}^{s}
(0,T).$$
Recalling the definition of operator $W_{bdr} (t)$ in Remark
\ref{Rem 2.1}, the solution $v$ of equation (\ref{2.2.1*}) can be
expressed by
\begin{equation*}
v(x,t) = W_{\mathbb{R}}(t) \bar{\phi}  - W_{bdr} (t)  \left(
(W_{\mathbb{R}}(t) \bar{\phi})  (0, t),   \; (W_{\mathbb{R}}(t)
\bar{\phi})_x (0, t) \right).
\end{equation*}
Thanks to  Proposition \ref{Pro 6.1} and estimate (\ref{2.1.2*}), we
obtain estimate (\ref{2.2.2*}).\endpf

\bigskip

\noindent {\bf Proof of estimate (\ref{2.3.2*}). }
By estimate (\ref{2.2.2*}), similar to Step 1 of the proof of
Proposition \ref{Pro 6.2}, we can obtain estimate (\ref{2.3.2*}).
\endpf


\subsection{Proof of estimates (\ref{2.1.3*})} \label{Sec 3.3}

\noindent {\bf Proof of estimate (\ref{2.1.3*}). }
By  (\ref{2.1.2*}), it suffices to prove that for $-2<s<0$, it holds
\begin{equation}\label{6.1.12}
 \|  t^{\frac{|s|} {4} + \varepsilon }  v \| _{ C([0,T];L^2 (\mathbb{R}^+)) }
+ \| t^{\frac{|s|} {4} + \varepsilon }  v\| _{L^2( 0,T; H^{2}
(\mathbb{R}^+ ) )}
\leq C  \left( \| \vec{h} \|_{ {\cal H}^s (0,T) } + \| t^{\frac{|s|}
{4} + \varepsilon } \vec{h} \| _{ {\cal H}^0 (0,T) }   \right).
\end{equation}

 It follows from (\ref{2.1.1*}) that $q :=
t^{\frac{|s|}{4} + \varepsilon} v$ satisfies
\begin{eqnarray}\label{6.1.13}
\left\{\begin{array}{ll}
 q_{t} + q_{xxxx} + \delta q_{xxx} = (\frac{|s|}{4} + \varepsilon)
 t^{\frac{|s|}{4} + \varepsilon -1 } v, & (x,t)\in \mathbb{R}^+
\times (0,T),
\\ \ns
q(x,0)=0, & x\in \mathbb{R}^+ ,
\\ \ns
q(0,t) = t ^{\frac{|s|}{4} + \varepsilon} h_1(t), \quad    q_x(0,t)
= t ^{\frac{|s|}{4} + \varepsilon} h_2(t), \qquad & t\in (0,T).
\end{array}\right.
\end{eqnarray}
Let $\theta$ and $\vartheta$ be solutions of
\begin{eqnarray}\label{6.1.15}
\left\{\begin{array}{ll}
\theta_{t} + \theta_{xxxx} + \delta \theta_{xxx} = (\frac{|s|}{4} +
\varepsilon) t^{\frac{|s|}{4} + \varepsilon -1 } v, \qquad &
(x,t)\in \mathbb{R}^+ \times (0,T),
\\ \ns
\theta(x,0)= 0, &  x\in \mathbb{R}^+,
\\ \ns
\theta(0,t) = 0, \quad    \theta_x(0,t) = 0,  & t\in (0,T),
\end{array}\right.
\end{eqnarray}
and
\begin{eqnarray}\label{6.1.17}
\left\{\begin{array}{ll} \displaystyle \vartheta_{t} +
\vartheta_{xxxx} + \delta \vartheta_{xxx} = 0, &   (x,t)\in
\mathbb{R}^+ \times (0,T),
\\ \ns
\vartheta(x,0)=0,  &  x\in \mathbb{R}^+ ,
\\ \ns
\vartheta(0,t) =   t^{\frac{|s|}{4} + \varepsilon } h_1(t),  \quad
\vartheta_x(0,t) = t^{\frac{|s|}{4} + \varepsilon } h_2(t), \qquad &
t\in (0,T).
\end{array}\right.
\end{eqnarray}
It follows from (\ref{6.1.13})-(\ref{6.1.17}) that
\begin{equation}\label{6.1.19}
q = \theta +\vartheta .
\end{equation}

From the results in  Chapter 4, Section 15 of \cite{Lions}, we get
that the solution of equation (\ref{2.1.1*}) satisfies $v \in
H^{\frac{1}{2} + \frac{s}{4}} (0,T; L^2(\mathbb{R}^+)) $ and
\begin{equation}\label{6.1.23}
\| v \| _{H^{\frac{1}{2} + \frac{s}{4}} (0,T; L^2(\mathbb{R}^+)) }
\displaystyle            \leq C    \| \vec{h} \|_{ {\cal H}^s (0,T)
}.
\end{equation}
For $-2<s <0$, Sobolev's embedding theorem shows that  $v \in
H^{\frac{1}{2} + \frac{s}{4}} (0,T; L^2(\mathbb{R}^+))
\hookrightarrow L^{\frac{4}{|s|}} (0,T; L^2(\mathbb{R}^+))$ and
\begin{equation}\label{6.1.25}
  \| v \| _{ L^{\frac{4}{|s|}} (0,T;
L^2(\mathbb{R}^+))}   \leq  C  \| v \| _{H^{\frac{1}{2} +
\frac{s}{4}} (0,T; L^2(\mathbb{R}^+)) } .
\end{equation}
Combining $t^{\frac{|s|}{4} + \varepsilon -1 } \in L^{\frac{4}{4 -
|s|}} (0,T)$ and  $v\in L^{\frac{4}{|s|}} (0,T; L^2(\mathbb{R}^+))$,
the H$\ddot{\text{o}}$lder's inequality shows that $(\frac{|s|}{4} +
\varepsilon)   t^{\frac{|s|}{4} + \varepsilon -1 } v \in L^1 (0,T;
L^2 (\mathbb{R}^+))$. That is
\begin{equation}\label{6.1.27}
\begin{array}{ll}
   \|(\frac{|s|}{4} + \varepsilon)    t^{\frac{|s|}{4} + \varepsilon -1 } v  \|  _{L^1 (0,T; L^2 (\mathbb{R}^+) )}
 \leq     C    \|  t^{\frac{|s|}{4} + \varepsilon -1 } v  \|  _{L^1 (0,T; L^2
(\mathbb{R}^+))} \\ \ns
 \leq C  \| t^{\frac{|s|}{4} + \varepsilon -1 }   \|_
{L^{\frac{4}{4 - |s|}} (0,T)}     \| v \| _{ L^{\frac{4}{|s|}} (0,T;
L^2(\mathbb{R}^+))}
 \leq C      \| v \| _{ L^{\frac{4}{|s|}} (0,T; L^2(0,L))}.
\end{array}
\end{equation}
Since $(\frac{|s|}{4} + \varepsilon)   t^{\frac{|s|}{4} +
\varepsilon -1 } v \in L^1 (0,T; L^2 (\mathbb{R}^+))$, applying
(\ref{2.3.2*}) in Proposition \ref{Pro 2.3*} to  equation
(\ref{6.1.15}),  we yield
\begin{equation}\label{6.1.33}
\begin{array}{ll}
\|\theta \|_{ C([0,T];L^2 (\mathbb{R}^+))}   +  \|\theta \|_{ L^2
(0,T; H^2 (\mathbb{R}^+))}  & \leq C    \|  (\frac{|s|}{4} +
\varepsilon) t^{\frac{|s|}{4} + \varepsilon -1 } v \|  _{L^1 (0,T;
L^2 (\mathbb{R}^+))}.
\end{array}
\end{equation}
Thus,  combining (\ref{6.1.23})-(\ref{6.1.33}), one has
\begin{equation}\label{6.1.35}
\begin{array}{ll}
\displaystyle  \|\theta \|_{ C([0,T];L^2 (\mathbb{R}^+))}   +
\|\theta \|_{ L^2 (0,T; H^2 (\mathbb{R}^+))}
  \leq C    \| \vec{h} \|_{ {\cal
H}^s (0,T) }.
\end{array}
\end{equation}

\smallskip

From $t^{\frac{|s|}{4} + \varepsilon} \vec{h} \in {\cal H} ^{0}
(0,T)$ and  estimate (\ref{2.1.2*}), we get that the solution of
equation (\ref{6.1.17}) satisfies
\begin{equation}\label{6.1.36}
\| \vartheta \|_{ C([0,T];L^2 (\mathbb{R}^+))}   +  \| \vartheta
\|_{ L^2 (0,T; H^2 (\mathbb{R}^+))}   \leq C   \| t^{\frac{|s|} {4}
+ \varepsilon } \vec{h} \| _{ {\cal H}^0 (0,T) } .
\end{equation}
Therefore, (\ref{6.1.19}),  (\ref{6.1.35}) and (\ref{6.1.36}) deduce
the desired
 estimate (\ref{6.1.12}).\endpf


\section{Well-posedness for $s\geq 0$} \label{Sec 4}

This section is devoted to proving Theorem \ref{Thm 1.2*}  for the
case $s\geq 0$. First of all, we deduce the local well-posedness
result of Kuramoto-Sivashinsky equation (\ref{1.1.1*}). Then,
combining  global priori estimate,  equation (\ref{1.1.1*}) is
globally well-posed.


\subsection{Local well-posedness for $s= 0$} \label{Sec 4.1}

The main purpose of this subsection is to discuss the local
well-posedness of Kuramoto-Sivashinsky equation (\ref{1.1.1*}).
By the result of the associated linear equation presented in
Subsection \ref{Sec 3.1}, if we further obtain a suitable estimate
for $u_{xx}+uu_x$ (see Lemma \ref{Lma 3.1} below), then from the
fixed point theory we can seek a fixed point solution.


\begin{lemma}\label{Lma 3.1}
There is a constant $C> 0$ such that for any  $ u$, $v \in  C (  [0,
T] ; L^2 (\mathbb{R}^+) ) \cap L^2 ( 0,T; H^{2} (\mathbb{R}^+)  )$,
it holds
\begin{equation*}
\begin{array}{ll}
\displaystyle  \big\| u_{xx}+ (u v)_x   \big\| _{L^1(0,T; L^2
(\mathbb{R}^+))}  \leq C \big[ T^{\frac{1}{2}}
  \|u\|  _{L^2 (0,T; H^{2} (\mathbb{R}^+) )}  +  \\ \ns
\displaystyle  \quad  ( T^{\frac{1}{2}} +  T^{\frac{1}{4}})
  \big( \| u \| _{C ( [0,T]; L^2 (\mathbb{R}^+)
) } +   \| u \| _{L^2 (0,T; H^{2} (\mathbb{R}^+) )} \big)   \big( \|
v \| _{C ( [0,T]; L^2 (\mathbb{R}^+) ) } +   \| v \| _{L^2 (0,T;
H^{2} (\mathbb{R}^+) )}  \big) \big].
\end{array}
\end{equation*}

\end{lemma}

{\it Proof.} On the one hand, by the H$\ddot{\text{o}}$lder's
inequality,
$$\int_0^T \|u_{xx} (\cdot, t) \| _{L^2 (\mathbb{R}^+)}  dt  \leq
\left( \int_0^T   dt\right)^{\frac{1}{2}}
\left( \int_0^T \|u_{xx} (\cdot, t) \|^2 _{L^2 (\mathbb{R}^+)}
dt\right)^{\frac{1}{2}}
\leq  C T^{\frac{1}{2}}
  \|u\|  _{L^2 (0,T; H^{2} (\mathbb{R}^+) )}.
$$

On the other hand, by the Gagliardo-Nirenberg's inequality, we get
\begin{equation}\label{3.1.8}
\begin{array}{ll}
\displaystyle  \int_0^T  \|  u(\cdot, t) v_x(\cdot, t) \|
_{L^2(\mathbb{R}^+)} dt
\displaystyle  \leq  \int_0^T   \|u(\cdot, t) \| _{L^{\infty}
(\mathbb{R}^+)} \|v_x (\cdot, t) \| _{L^2 (\mathbb{R}^+)} dt \\ \ns
\displaystyle  \leq  C  \int _0^T   \left(   \| u(\cdot, t) \|
_{L^2(\mathbb{R}^+)}
  +   \| u(\cdot, t) \| _{L^2(\mathbb{R}^+)} ^{\frac{1}{2}}   \| u_x(\cdot, t) \| _{L^2(\mathbb{R}^+)}
  ^{\frac{1}{2}}\right)  \| v_x(\cdot, t) \| _{L^2(\mathbb{R}^+)}  dt.
\end{array}
\end{equation}
From the H$\ddot{\text{o}}$lder's inequality, (\ref{3.1.8}) gives
\begin{equation*}
\begin{array}{ll}
\displaystyle  \int_0^T  \|  u(\cdot, t) v_x(\cdot, t) \|
_{L^2(\mathbb{R}^+)} dt\\ \ns
\leq \displaystyle C \sup_{t\in [0,T]} \|  u(\cdot, t)
\|_{L^2(\mathbb{R}^+)}
  \left( \int_0^T dt
 \right)^{\frac{1}{2}}
   \left( \int_0^T    \| v_x(\cdot, t) \| _{L^2(\mathbb{R}^+)} ^2 dt  \right)
  ^{\frac{1}{2}}   \\  \ns
\displaystyle \quad + C \sup_{t\in [0,T]} \|  u(\cdot, t)
\|_{L^2(\mathbb{R}^+)} ^{\frac{1}{2}}    \left( \int_0^T
dt\right)^{\frac{1}{4}} \left( \int_0^T \| u_x(\cdot,t) \| _{L^2
(\mathbb{R}^+)} ^{2} dt
 \right)^{\frac{1}{4}}
   \left( \int_0^T    \| v_x(\cdot, t) \| _{L^2(\mathbb{R}^+)} ^{2} dt  \right)
  ^{\frac{1}{2}}   \\  \ns
\displaystyle  \leq C    T^{\frac{1}{2}} \| u \| _{C ( [0,T]; L^2
(\mathbb{R}^+) ) } \| v \| _{L^2 (0,T; H^{2} (\mathbb{R}^+) )} \\
\ns
\displaystyle \quad   + C  T^{\frac{1}{4}}    \| u
\|^{\frac{1}{2}}  _{C ( [0,T]; L^2 (\mathbb{R}^+) )
  } \| u \|^{\frac{1}{2}}  _{L^2 (0,T; H^{2} (\mathbb{R}^+) )}
\| v \| _{L^2 (0,T; H^{2} (\mathbb{R}^+) )}.
\end{array}
\end{equation*}
This completes the proof of Lemma \ref{Lma 3.1}.\endpf

\begin{proposition}\label{Pro 3.1}
Let $T>0$, $\phi \in L^2 (\mathbb{R}^+) $ and $(h_1, h_2) \in
H^{\frac{3}{8}} (0,T) \times H^{\frac{1}{8}} (0,T)$, then there
exists a $T_* \in (0,T]$ depending on $ \|\phi \|_{
L^2(\mathbb{R}^+)} +   \|h_1 \|_{ H^{\frac{3}{8}} (0,T) }  +   \|
h_2 \|_{ H^{\frac{1}{8}} (0,T) } $, such that equation
(\ref{1.1.1*}) admits a unique solution $u\in C ([0, T_* ];
L^2(\mathbb{R}^+))\cap L^2 ( 0,T_* ; H^{2} (\mathbb{R}^+) )$ with
$ \partial_x^k u \in C ( [0, +\infty); H^{\frac{3}{8} - \frac{k}{4}}
(0, T_* ) )$ ($k=0,1$).
Moreover, the corresponding solution map from the space of initial
and boundary data to the solution space is  continuous.

\end{proposition}

{\it Proof. } The solution of equation (\ref{1.1.1*}) can be written
in the form
$$u(t) = W_c(t) \phi + W_{bdr} (t) \vec{h}  - \int_0^t W_c(t-\tau) \left( u_{xx} + u u_x  \right) (\tau)d\tau$$
(recall the definitions of $W_{bdr} (t)$ and $W_c(t)$ in Remark
\ref{Rem 2.1} and Remark \ref{Rem 2.2}).  Define
\begin{equation}\label{3.1.10}
\Gamma (w) = W_c(t) \phi + W_{bdr} (t) \vec{h}  - \int_0^t
W_c(t-\tau) \left( w_{xx} + w w_x  \right) (\tau)d\tau.
\end{equation}
By Proposition \ref{Pro 2.1*} - \ref{Pro 2.3*} and Lemma \ref{Lma
3.1}, for any $w\in X_{0, T_*} (d) : = \left \{w\in X_{0, T_*} \;
\big|\ \|w\| _{X_{0, T_*} } \leq d \right \}$,
\begin{equation}\label{3.1.17}
\displaystyle
  \|\Gamma (w) \|_{X_{0, T_* } }
\leq   \displaystyle  C_1  \left(   \|\phi \|_{ L^2(\mathbb{R}^+)} +
\| \vec{h} \|_{ {\cal H}^{0} (0,T) }   \right)
 \displaystyle  +  C_2  \left[ T_*^{\frac{1}{2}} \| w \| _{X_{0, T_*
} }  +  ( T_*^{\frac{1}{2}} + T_*^{\frac{1}{4}} ) \| w \| ^2_{X_{0,
T_* } } \right].
\end{equation}
Set
\begin{equation}\label{3.1.19}
d= 2 C_1   \left(   \|\phi \|_{ L^2(\mathbb{R}^+)} +   \| \vec{h}
\|_{ {\cal H}^{0} (0,T) } \right)
\end{equation}
and Choose $T_*$ such that
\begin{equation}\label{3.1.20}
\left\{\begin{array}{ll} C_2  T_* ^{\frac{1}{2}} \leq \frac{1}{4},
\\ \ns
C_2   ( T_* ^{\frac{1}{2}} + T_* ^{\frac{1}{4}}  ) d \leq
\frac{1}{4}.
\end{array}\right.
\end{equation}
By (\ref{3.1.17})-(\ref{3.1.20}), we have that for any
$\displaystyle w \in  X_{0, T_*} (d) $,
$$ \|\Gamma (w) \|_{X_{0, T_* } } \leq d,$$
which means that $\Gamma$ maps $X_{0, T_* } (d)$ into $X_{0, T_* }
(d)$.

\smallskip

Furthermore, $\Gamma$ is a contraction map of $X_{0, T_* } (d)$. In
fact, from (\ref{3.1.10}), for any $w_1$, $ w_2\in X_{0, T_* } (d)$,
\begin{equation}\label{3.1.24}  \Gamma(w_1)   -   \Gamma(w_2)   =
- \int_0^t   W_c(t-\tau) \Big[ (w_1 - w_2 ) _{xx} + \frac{1}{2}
\left[(w_1  + w_2 )  ( w_1 - w_2) \right] _x \Big] (\tau)d\tau.
\end{equation}
By Proposition \ref{Pro 2.3*} and Lemma \ref{Lma 3.1},
(\ref{3.1.24}) leads to
\begin{equation}\label{3.1.27}
\begin{array}{ll}
\displaystyle   \|  \Gamma(w_1)   -   \Gamma(w_2)  \|_{X_{0, T_* } } \\
\ns
\displaystyle   \leq   C_2  \left[  T_* ^{\frac{1}{2}}
  \|w_1 - w_2\|  _{X_{0,T_* }}   + \frac{1}{2}    (
T_*^{\frac{1}{2}} +    T_*^{\frac{1}{4}}  )
   \|w_1 + w_2 \|  _{X_{0, T_* }} \| w_1 - w_2 \|  _{X_{0,  T_* }}  \right].
\end{array}
\end{equation}
Then (\ref{3.1.20}) and (\ref{3.1.27}) imply that for any $w \in
X_{0, T_* } (d)$,
$$ \|  \Gamma(w_1)   -   \Gamma(w_2)  \|_{X_{0, T_* } }
\leq   \frac{1}{2}
  \|w_1 - w_2 \|  _{X_{0,  T_* }} .$$

Since $\Gamma$ is a contraction mapping from $X_{0, T_* } (d)$ to
$X_{0, T_* } (d)$, the fixed point theorem shows that $\Gamma$
exists a fixed point $u$. Then equation (\ref{1.1.1*}) admits a
unique solution $u=\Gamma(u)$ in $X_{0, T_* } (d)$.
This completes the proof of Proposition \ref{Pro 3.1}.\endpf

\subsection{Global well-posedness for $s\geq 0$} \label{Sec 4}

This subsection is addressed to proving Theorem \ref{Thm 1.2*} for
$s\geq 0$. We divide the proof into three steps. In Step 1 and Step
2, Theorem \ref{Thm 1.2*} is shown to be true for $s=0$ and $s=4$,
respectively. Then by nonlinear interpolation theory, we prove it is
true for any $s\in(0,4)$ in Step 3.
If $s\in (4,8)$, then $s-4\in (0,4)$. Thus, by applying the result
of Step 3, Theorem \ref{Thm 1.2*} holds for any $s\in (4,8)$. The
same method can be used for $s\geq 8$.

To state our results clearly and concisely, we introduce the
following notations.
\begin{equation}\label{1.1.2*}
\begin{array}{ll}
\displaystyle Z _{s,T}  = H^s (\mathbb{R}^+) \times  {\cal H} ^s
(0,T), \qquad \| (\phi, \vec{h}) \|_{Z_{s,T} } =  \| \phi
\|_{H^s(\mathbb{R}^+)} + \| \vec{h} \| _{{\cal H} ^s (0,T)} .
\end{array}
\end{equation}

\medskip

{\it Step 1.} By the local well-posedness results, to prove Theorem
\ref{Thm 1.2*} for the case $s= 0$, we only need the following
global priori estimate (\ref{4.3}) for smooth solutions of equation
(\ref{1.1.1*}). Then we can patch local solutions together to a
global solution.

For given $T>0$, we claim that there exists a continuous
nondecreasing function $\gamma : \mathbb{R}^+   \rightarrow
\mathbb{R}^+$ such that for any smooth solution $u$ of equation
(\ref{1.1.1*}), it hold
\begin{equation}\label{4.3}
 \| u \|_{C( [0,T]; L^2(\mathbb{R}^+ ) ) }   +  \| u \|_{ L^2
(  0,T ; H^{2} (\mathbb{R}^+ )  )}  \leq  \gamma \left( \| (\phi,
\vec{h}) \|_{Z_{0,T} }  \right).
\end{equation}
In fact, let $u$ be a smooth solution of equation (\ref{1.1.1*}),
let $y$ and $z$ solve
\begin{eqnarray}\label{4.1.5}
\left\{\begin{array}{ll} \displaystyle y_{t} + y_{xxxx} + \delta
y_{xxx} = 0,  & (x,t)\in \mathbb{R}^+  \times \mathbb{R}^+,
\\ \ns
y(x,0)= 0 ,  &  x\in \mathbb{R}^+,
\\ \ns
y(0,t)=h_1(t),\quad   y_x(0,t)=h_2(t), \qquad  &  t\in \mathbb{R}^+,
\end{array}\right.
\end{eqnarray}
and
\begin{eqnarray}\label{4.1.7}
\left\{\begin{array}{ll} \displaystyle z_{t} + z_{xxxx} + \delta
z_{xxx} + z_{xx} + z z_x = -(z y)_x - yy_x - y_{xx}, \qquad &
(x,t)\in \mathbb{R}^+ \times \mathbb{R}^+,
\\ \ns
z(x,0)= \phi (x)  ,  &  x\in \mathbb{R}^+,
\\ \ns
z(0,t)=0,\quad   z_x(0,t)=0,  &   t\in \mathbb{R}^+.
\end{array}\right.
\end{eqnarray}
Then it is easy to check that $u= y+z$.

\medskip

Multiplying both side of the first equation in (\ref{4.1.7}) by $z$
and integrating by parts, we get
\begin{equation}\label{4.1.13}
\begin{array}{ll}
\displaystyle   \frac{1}{2} \frac{d}{dt}   \|  z \| ^2
_{L^2(\mathbb{R}^+)} + \| z_{xx} \|^2 _{L^2(\mathbb{R}^+)}  \\ \ns
\displaystyle    =   \| z_{x} \|^2 _{L^2(\mathbb{R}^+)}  + \int_0^{+
\infty} z_x y z dx
 + \frac{1}{2} \int_0^{+ \infty}  z_x  y ^2  dx    - \int_0^{+ \infty}   z_{xx} y dx.
\end{array}
\end{equation}
Now we estimate the four terms on the right hand side of
(\ref{4.1.13}). The Gagliardo - Nirenberg's inequality and the
Young's inequality imply that for any $\varepsilon>0$,
\begin{equation}\label{4.1.14}
\begin{array}{ll}
\displaystyle \| z_x \| ^2 _{L^2(\mathbb{R}^+)}     \leq C_1 \left(
\|z_{xx} \| ^{\frac{1}{2}}_{L^2(\mathbb{R}^+)} \|z \|
^{\frac{1}{2}}_{L^2(\mathbb{R}^+)}  +  \|z \| _{L^2(\mathbb{R}^+)}
\right) ^2 \\ \ns
 \displaystyle  \leq  2C_1  \|z_{xx} \| _{L^2(\mathbb{R}^+)}   \|z
\| _{L^2(\mathbb{R}^+)}   + 2C_1  \|z \| ^2_{L^2(\mathbb{R}^+)}\\
\ns
  \leq  \varepsilon^2 \|z_{xx} \| ^2_{L^2(\mathbb{R}^+)} +   \left(
\frac{C_1^2}{\varepsilon^2} + 2C_1 \right) \|z \|
^2_{L^2(\mathbb{R}^+)}.
\end{array}
\end{equation}
From  the Gagliardo - Nirenberg's inequality, the
H$\ddot{\text{o}}$lder's inequality, the Young's inequality and
(\ref{4.1.14}), we get
\begin{equation}\label{4.1.16}
\begin{array}{ll}
\displaystyle  \int_0^{+ \infty}  z_x  y z dx \leq \|z_x
\|_{L^{\infty}(\mathbb{R}^+)} \| y z \|_{L^1(\mathbb{R}^+)}
\displaystyle   \leq  \|z_x \|_{L^{\infty}(\mathbb{R}^+)} ^2   + \|
y z \|_{L^1(\mathbb{R}^+)} ^2  \\ \ns
\displaystyle   \leq C_2 \left( \| z_{xx} \| _{L^2(\mathbb{R}^+)}
^{\frac{1}{2}}   \| z_{x} \| _{L^2(\mathbb{R}^+)}
  ^{\frac{1}{2}}  +   \| z _x \| _{L^2(\mathbb{R}^+)}  \right)  ^2     +    \| y  \|^2 _{L^2(\mathbb{R}^+)}    \| z \|^2_{L^2(\mathbb{R}^+)}   \\ \ns
\displaystyle   \leq 2 C_2   \| z_{xx} \| _{L^2(\mathbb{R}^+)}   \| z_{x} \| _{L^2(\mathbb{R}^+)}
 +    2 C_2  \| z _x \| ^2_{L^2(\mathbb{R}^+)}    +    \| y  \|^2 _{L^2(\mathbb{R}^+)}    \| z \|^2_{L^2(\mathbb{R}^+)}   \\  \ns
  \leq  \varepsilon \|z_{xx} \| ^2_{L^2(\mathbb{R}^+)} +
\left(  \frac{C_2^2} {\varepsilon }  +  2C_2  \right) \|z_x \|
^2_{L^2(\mathbb{R}^+)}  +    \| y  \|^2 _{L^2(\mathbb{R}^+)} \| z
\|^2_{L^2(\mathbb{R}^+)}    \\  \ns
  \leq  \left( \varepsilon + \varepsilon  C_2^2  +2 \varepsilon^2  C_2  \right)  \|z_{xx} \| ^2_{L^2(\mathbb{R}^+)} +
 \left(  \frac{C_1^2}{\varepsilon^2} + 2C_1 \right)  \left(  \frac{C_2^2} {\varepsilon }  +  2C_2  \right)  \|z \| ^2_{L^2(\mathbb{R}^+)}  \\ \ns
\quad\displaystyle +    \| y  \| ^2_{L^2(\mathbb{R}^+)} \| z
\|^2_{L^2(\mathbb{R}^+)}.
\end{array}
\end{equation}
Similar to (\ref{4.1.16}), one yields
\begin{equation}\label{4.1.18}
\begin{array}{ll}
\displaystyle   \frac{1}{2} \int_0^{+ \infty}  z_x  y ^2  dx
 \leq \|z_x \|_{L^{\infty}(\mathbb{R}^+)}
\| y \|^2_{L^2(\mathbb{R}^+)}    \leq   \|z_x
\|_{L^{\infty}(\mathbb{R}^+)} ^2   + \| y \|_{L^2(\mathbb{R}^+)} ^4
\\ \ns
\leq  \left( \varepsilon + \varepsilon  C_2^2  +2 \varepsilon^2  C_2
\right)  \|z_{xx} \| ^2_{L^2(\mathbb{R}^+)}  +
 \left(  \frac{C_1^2}{\varepsilon^2} + 2C_1 \right)  \left(  \frac{C_2^2}
{\varepsilon } + 2C_2 \right)  \|z \| ^2_{L^2(\mathbb{R}^+)}   \\
\ns
 \quad + \| y \|_{L^2(\mathbb{R}^+)} ^4.
\end{array}
\end{equation}
Again, by  the Young's inequality, it holds
\begin{equation}\label{4.1.20}
 - \int_0^{+ \infty}    z_{xx} y dx   \leq   \|z_{xx}\| _{L^2(\mathbb{R}^+)}  \|y\| _{L^2(\mathbb{R}^+)}     \leq  \varepsilon \|z_{xx} \|
^2_{L^2(\mathbb{R}^+)} + \frac{1}{4 \varepsilon} \|y \|
^2_{L^2(\mathbb{R}^+)}.
\end{equation}
Choosing $\varepsilon$ small enough, from (\ref{4.1.13}) -
(\ref{4.1.20}), we deduce
\begin{equation}\label{4.1.21}
 \frac{d}{dt}   \| z  \|^2 _{L^2(\mathbb{R}^+)}  +   \| z _{xx} \|^2 _{L^2(\mathbb{R}^+)}
\displaystyle    \leq  \left( C  +   \| y \|^2 _{L^2(\mathbb{R}^+)}
\right) \| z \|^2 _{L^2(\mathbb{R}^+)}
 +   \| y \|^4 _{L^2(\mathbb{R}^+)}  +C  \| y \|^2 _{L^2(\mathbb{R}^+)}.
\end{equation}

Proposition \ref{Pro 2.1*} implies that
\begin{equation}\label{4.1.23}
\| y \| _{ C([0,T]; L^2(\mathbb{R}^+) )}   +   \| y \|_{ L^2  (  0,T
; H^{2} (\mathbb{R}^+)  )}     \leq  C \| \vec{h} \|_{{\cal H} ^0
(0,T)}.
\end{equation}
From (\ref{4.1.21}) and (\ref{4.1.23}), for any $t\in [0,T]$,
\begin{equation}\label{4.1.24}
\begin{array}{ll}
\displaystyle \frac{d}{dt}   \| z (\cdot,t) \| ^2_{L^2(\mathbb{R}^+)}   +   \| z_{xx} (\cdot,t) \| ^2_{L^2(\mathbb{R}^+)}\\
  \ns
\displaystyle    \leq  C  \left( 1 +
 \| \vec{h}\|^2_{{\cal H} ^0 (0,T)} \right) \| z (\cdot, t) \|^2_{L^2(\mathbb{R}^+)}     +    C  \left(   \| \vec{h}\|^4_{{\cal H} ^0 (0,T)} +
 \| \vec{h}\|^2_{{\cal H} ^0 (0,T)} \right) .
\end{array}
\end{equation}
Applying the Gronwall's inequality to (\ref{4.1.24}), one gets
%
%
%
\begin{equation}\label{4.1.27}
\begin{array}{ll}
\displaystyle  \sup_{t\in [0,T]} \| z (\cdot,t) \| ^2_{L^2(\mathbb{R}^+)}  \\
\ns
\displaystyle \leq e^{CT \left( 1 +
 \| \vec{h}\|^2_{{\cal H} ^0 (0,T)} \right)}  \left[ \| \phi\| ^2_{L^2(\mathbb{R}^+)}  +  CT \left(   \| \vec{h}\|^4_{{\cal H} ^0 (0,T)} +
 \| \vec{h}\|^2_{{\cal H} ^0 (0,T)} \right)   \right].
\end{array}
\end{equation}
Integrating (\ref{4.1.24}) on $[0,T]$, we get
\begin{equation*}
\begin{array}{ll}
\displaystyle   \|z_{xx} \|^2 _{L^2(0,T; L^2(\mathbb{R}^+ ))} \\
\ns
\displaystyle  \leq    \| \phi\| ^2_{L^2(0,L)} +    C  \left( 1 +
 \| \vec{h}\|^2_{{\cal H} ^0 (0,T)} \right)\int_0^T \|z(\cdot, t) \|^2_{L^2 (\mathbb{R}^+ )}  dt
  +    C T  \left(   \| \vec{h}\|^4_{{\cal H} ^0 (0,T)} +
 \| \vec{h}\|^2_{{\cal H} ^0 (0,T)} \right)   \\ \ns
\displaystyle  \leq    \| \phi\| ^2_{L^2(0,L)} +    C  T \left( 1 +
 \| \vec{h}\|^2_{{\cal H} ^0 (0,T)} \right)   \sup_{t\in [0,T]} \| z (\cdot,t) \| ^2_{L^2(\mathbb{R}^+)}
  +    C T  \left(   \| \vec{h}\|^4_{{\cal H} ^0 (0,T)} +
 \| \vec{h}\|^2_{{\cal H} ^0 (0,T)} \right) ,
\end{array}
\end{equation*}
which combining (\ref{4.1.27}) gives that there exists a continuous
nondecreasing function $\gamma _* : \mathbb{R}^+   \rightarrow
\mathbb{R}^+$ such that
\begin{equation}\label{4.1.28}
\| z \|_{ L^2  (  0,T ; H^{2} (\mathbb{R}^+)  )} \leq \gamma _*
\left( \| (\phi, \vec{h}) \|_{Z_{0,T} } \right).
\end{equation}
From (\ref{4.1.23}), (\ref{4.1.27}) and (\ref{4.1.28}) and $u= y+z$,
it easy to see that (\ref{4.3}) holds.

\medskip

{\it Step 2.} We now prove that Theorem \ref{Thm 1.2*} is true for
$s=4$. By the result of Step 1 and $ (\phi, \vec{h}) \in  Z_{4,T}
\subset  Z_{0,T}$, the solution of equation (\ref{1.1.1*}) satisfies
(recall the definition of $X_{0,T}$ in (\ref{2.3*}))
\begin{equation}\label{4.2.3}
 \| u \|_{X_{0,T}} \leq  C \| (\phi,
\vec{h}) \|_{Z_{0,T} }.
\end{equation}
Furthermore, it is easy to check that $y = u_t$ solves
\begin{eqnarray}\label{4.2.5}
\left\{\begin{array}{ll} \displaystyle y_{t} + y_{xxxx} + \delta
y_{xxx} = - y_{xx} - (uy)_x ,  &   (x,t)\in  \mathbb{R}^+ \times
\mathbb{R}^+,
\\ \ns
y(x,0)= \psi (x) = - \phi ^{\prime\prime\prime\prime} (x) - \delta
\phi ^{\prime\prime\prime} (x) - \phi ^{\prime\prime} (x)  - \phi
(x) \phi ^{\prime} (x),  \qquad & x\in \mathbb{R}^+,
\\ \ns
y(0,t)=h_1 ^{\prime} (t),\quad   y_x(0,t)=h_2 ^{\prime}  (t),  &
 t\in \mathbb{R}^+.
\end{array}\right.
\end{eqnarray}
By Proposition \ref{Pro 2.1*} - \ref{Pro 2.3*} and Lemma \ref{Lma
3.1}, there is a constant $C>0$ such that for any $T_*\in (0,T]$,
\begin{equation}\label{4.2.8}
\begin{array}{ll}
\| y \| _{ X_{0,T_*}} \leq   C   \| (\psi, \vec{h}^{\prime})
\|_{Z_{0,T} } +    C  \left[ T_* ^{\frac{1}{2}} \| y \| _{X_{0,T_*}
}  +  ( T_* ^{\frac{1}{2}} + T_*^{\frac{1}{4}} ) \| u \| _{X_{0,T} }
\| y \| _{X_{0,T_*} } \right].
\end{array}
\end{equation}
Choose $T_*$ such that
\begin{equation}\label{4.2.10}
\left\{\begin{array}{ll}  C  T_*^{\frac{1}{2}} \leq \frac{1}{4},
\\ \ns
C  ( T_*^{\frac{1}{2}} + T_*^{\frac{1}{4}} ) \| u \| _{X_{0,T}} \leq
\frac{1}{4}.
\end{array}\right.
\end{equation}
Combining (\ref{4.2.8}) and (\ref{4.2.10}), one has
\begin{equation}\label{4.2.12}
\| y \| _{ X_{0,T_*}} \leq  2 C  \| (\psi, \vec{h}^{\prime})
\|_{Z_{0,T} }.
\end{equation}
By (\ref{4.2.10}),  $T_*$ only depends on $\| u \| _{X_{0,T}}$.
Furthermore, combining (\ref{4.2.3}), we see that $T_*$ only depends
on $\| (\phi, \vec{h}) \|_{Z_{0,T} } $.
Recalling (\ref{4.2.5}), by a standard density argument,
(\ref{4.2.12}) shows that there exists a continuous nondecreasing
function $C_0 : \mathbb{R}^+ \rightarrow \mathbb{R}^+ $ such that
\begin{equation*}
\begin{array}{ll}
\| y \| _{ X_{0,T}}  \leq
 C _0 \left( \| (\phi,
\vec{h}) \|_{Z_{0,T} } \right)     \| (\phi, \vec{h}) \|_{Z_{4,T}
 } .
\end{array}
\end{equation*}
Then Theorem \ref{Thm 1.2*} holds for $s=4$ since $y=-u_{xxxx} -
\delta u_{xxx} - u_{xx} - u u_x.$

\medskip

{\it Step 3.} Since Theorem \ref{Thm 1.2*} hold for $s=0$ and $s=4$,
we now use nonlinear interpolation theory (c.f. Theorem 1 in
\cite{Bona} or Theorem 2 in \cite{Tartar}) to prove it  for any $s
\in (0,4)$. Set
$$A_1 = Z_{0,T}, \quad  A_0 = Z_{4, T}, \quad  B_1 = C([0,T]; L^2(\mathbb{R}^+)), \quad  B_0 = C([0,T]; H^{4}(\mathbb{R}^+)).$$
%
%
%
Then for given $s\in (0,4)$,
\begin{equation}\label{4.3.2} [A_0, A_1] _{\frac{s}{4},
2} = Z_{s,T}, \qquad [B_0, B_1] _{\frac{s}{4}, 2} = C([0,T];
H^{s}(\mathbb{R}^+)).
\end{equation}

Let $\cal T$ be the solution map of equation (\ref{1.1.1*}), i.e. $u
= {\cal T} (\phi, \vec{h})$.
By the result proved in Step 2, we see that
 ${\cal T}: Z_{4,T}   \rightarrow
C([0,T]; H^{4}(\mathbb{R}^+)) $ and for any $(\phi, \vec{h}) \in
Z_{4,T}$,
\begin{equation}\label{4.3.3}
\| u \| _{ C([0,T]; H^{4}(\mathbb{R}^+)) }   \leq  C_1 \left( \|
(\phi, \vec{h}) \| _{Z_{0,T}} \right) \| (\phi, \vec{h}) \|
_{Z_{4,T}},
\end{equation}
where $C_1 : \mathbb{R}^+ \rightarrow \mathbb{R}^+$ is a continuous
nondecreasing function.

\medskip

Denote $u_1 = {\cal T} (\phi_1, \vec{h}_1)$, $u_2 = {\cal T}
(\phi_2, \vec{h}_2)$, $z=u_1 -u_2$ and $w =\frac{1}{2} (u_1 +u_2) $.
Then $z$ solves
\begin{eqnarray}\label{4.3.5}
\left\{\begin{array}{ll} \displaystyle z_{t} + z_{xxxx} + \delta
z_{xxx} + z_{xx} = - (w z)_x ,  & (x,t)\in  \mathbb{R}^+\times
\mathbb{R}^+,
\\ \ns
z(x,0)=\phi_1(x) - \phi _2 (x)  ,  &   x\in \mathbb{R}^+,
\\ \ns
\left( z(0,t), \; z_x(0,t) \right) = \vec{h}_{1} (t) - \vec{h}_{2 }
(t) ,  \qquad & t\in \mathbb{R}^+,
\end{array}\right.
\end{eqnarray}
By Proposition \ref{Pro 2.1*} - \ref{Pro 2.3*} and Lemma \ref{Lma
3.1}, there is a constant $C>0$ such that for any $T_*\in (0,T]$,
\begin{equation*}
\| z \| _{ X_{0,T^*}} \leq
  \displaystyle  C  \|
(\phi_1, \vec{h} _1) - (\phi_2, \vec{h} _2)\| _{Z_{0,T}}
  +     C    \left[
T_*^{\frac{1}{2}} \| z \| _{X_{0,T_*} } + ( T_*^{\frac{1}{2}} +
T_*^{\frac{1}{4}} ) \| w\| _{X_{0,T} } \| z \| _{X_{0,T_*} }
\right].
\end{equation*}
Choose $T_* $ such that
\begin{equation}\label{4.3.10}
\left\{\begin{array}{ll}  C   T_*^{\frac{1}{2}} \leq \frac{1}{4},
\\ \ns
C  ( T_*^{\frac{1}{2}} + T_*^{\frac{1}{4}}  ) \| w \| _{X_{0,T}}
\leq \frac{1}{4}.
\end{array}\right.
\end{equation}
Therefore,
\begin{equation}\label{4.3.12}
\| z \| _{ X_{0,T^*}} \leq 2 C   \| (\phi_1, \vec{h} _1)  - (\phi_2,
\vec{h} _2)\| _{Z_{0,T}}.
\end{equation}
By (\ref{4.3.10}), we know that  $T_*$ only depends on $\| w \|
_{X_{0,T}}$. Since $w=\frac{1}{2} (u_1 +u_2)$ and Theorem \ref{Thm
1.2*} is true for $s=0$, we see that $T^*$ only depends on $\|
(\phi_1 , \vec{h}_1) \| _{Z_{0,T}} +
 \| (\phi_2 , \vec{h}_2) \| _{Z_{0,T}}$.
Recalling (\ref{4.3.5}), (\ref{4.3.12}) implies that there exists a
continuous nondecreasing function $C_2 : \mathbb{R}^+ \rightarrow
\mathbb{R}^+ $ such that
\begin{equation*}
\| z \| _{ Y_{0,T}} \leq  C_2 \left( \| (\phi_1, \vec{h} _1) \|
_{Z_{0,T}}  + \| (\phi_2, \vec{h} _2)\| _{Z_{0,T}} \right) \|
(\phi_1, \vec{h} _1)  - (\phi_2, \vec{h} _2)\| _{Z_{0,T}}.
\end{equation*}
Hence,  we see that
 ${\cal T}: Z_{0,T}   \rightarrow
C([0,T]; L^2 (\mathbb{R}^+))  $ and for any $(\phi_1, \vec{h}_1)$,
$(\phi_2, \vec{h}_2) \in Z_{0,T}$,
\begin{equation}\label{4.3.15}
\displaystyle   \| u_1 - u_2 \| _{ C([0,T]; L^2 (\mathbb{R}^+)) }
\displaystyle   \leq  C_2 \left( \| (\phi_1, \vec{h} _1) \|
_{Z_{0,T}}  + \| (\phi_2, \vec{h} _2)\| _{Z_{0,T}} \right) \|
(\phi_1, \vec{h} _1)  - (\phi_2, \vec{h} _2)\| _{Z_{0,T}}.
\end{equation}
By (\ref{4.3.2}), (\ref{4.3.3}), (\ref{4.3.15}), the nonlinear
interpolation theory shows that Theorem \ref{Thm 1.2*} is true for
$s \in (0,4)$.
\endpf


\section{Well-posedness  for $-2 < s< 0$}

This section is devoted to proving Theorem \ref{Thm 1.2*}  for the
case $-2 < s< 0$. First of all, we deduce a suitable estimate for
$u_{xx}+uu_x$ (see Lemma \ref{Lma 5.1} below). Then we prove the
local well-posedness of Kuramoto-Sivashinsky equation (\ref{1.1.1*})
(see Proposition \ref{Pro 6.5}). Finally, we show that  equation
(\ref{1.1.1*}) admits a unique global solution.


\begin{lemma}\label{Lma 5.1}
Let $T>0$, $\varepsilon>0$ and $-2 + 4 \varepsilon < s< 0$. Then
there exists a constant $C> 0$ such that for any $ u$, $v \in
X_{s,T}^{\varepsilon}$ (recall the definition of
$X_{s,T}^{\varepsilon}$ in (\ref{6.1.1})), it holds
\begin{equation}\label{6.2.4}
\left \| u_{xx} \right\| _{L^1(0,T; H ^{s }(\mathbb{R} ^+ ))}   +
\|t^{ \frac{|s|}{4} + \varepsilon} u_{xx} \| _{L^1(0,T; L^2
(\mathbb{R} ^+ ))} \leq  C  T^{\alpha(s,\varepsilon)}
  \|u\|  _{X_{s,T}^{\varepsilon}}
\end{equation}
and
\begin{equation}\label{6.2.5}
\| u v_x \| _{L^1(0,T; H ^{s} (\mathbb{R} ^+ ))}    +   \| t^{
\frac{|s|}{4} + \varepsilon}  u v_x   \| _{L^1(0,T; L^2
(\mathbb{R}^+ ))} \leq C T^{\alpha(s,\varepsilon)}
   \|u\|  _{X_{s,T}^{\varepsilon}} \|v\|  _{X_{s,T}^{\varepsilon}} ,
\end{equation}
where $\alpha(s, \varepsilon): (-2,0)\times (0,+\infty) \rightarrow
(0,+\infty)$ is a continuous
 function.

\end{lemma}

{\it Proof.} We divide the proof into three steps. In Step 1,
estimate (\ref{6.2.4}) is shown to be true. In Step 2, estimate
$\| u v _x \| _{L^1(0,T; H ^{s} (\mathbb{R} ^+ ))}    \leq C
T^{\alpha(s, \varepsilon)}
  \|u\|  _{X_{s,T}^{\varepsilon}} \|v\|  _{X_{s,T}^{\varepsilon}} $
is proven for cases $-1 \leq s <0$, $- \frac{3}{2} <s <-1$, $s = -
\frac{3}{2}$ and $-2 +4 \varepsilon <s < - \frac{3}{2} $,
respectively.
Finally, we prove
$ \| t^{\frac{|s|}{4} + \varepsilon}  u v_x   \| _{L^1(0,T; L^2
(\mathbb{R} ^+ ))} \leq C T^{\alpha(s, \varepsilon)}
   \|u\|  _{X_{s,T}^{\varepsilon}} \|v\|  _{X_{s,T}^{\varepsilon}} $
in Step 3.

\smallskip

{\it Step 1.} We claim that (\ref{6.2.4}) holds.
Indeed, by (\ref{6.1.1}) and the H$\ddot{\text{o}}$lder's
inequality,
$$\int_0^T \|  t^{\frac{|s|}{4} + \varepsilon}  u_{xx}  \| _{L^2(\mathbb{R}^+ )}  dt
\leq  \left( \int_0^T   dt\right)^{\frac{1}{2}}
\left( \int_0^T   \| t^{ \frac{|s|}{4} +  \varepsilon}   u_{xx} \|^2
_{L^2(\mathbb{R} ^+ )} dt\right)^{\frac{1}{2}}
\leq   T^{\frac{1}{2}}
  \|u\|  _{X_{s,T}^{\varepsilon}}.
$$

From the H$\ddot{\text{o}}$lder's inequality  and the fact that the
embedding $L^2 (\mathbb{R}^+) \subset H^{s} (\mathbb{R}^+)$  is
continuous for $s<0$, we have that for  $ \frac{|s|}{2} + 2
\varepsilon<1$, i.e. $- 2 + 4 \varepsilon<s < 0$, it holds
\begin{equation*}
\begin{array}{ll}
\displaystyle   \int_0^T \|u_{xx}  \| _{H^{s }(\mathbb{R} ^+ )} dt
\leq
\left(  \int_0^T  t^{- \frac{|s|}{2} - 2 \varepsilon}
dt\right)^{\frac{1}{2}}
\left( \int_0^T  \| t^{ \frac{|s|}{4} +  \varepsilon} u_{xx}  \|^2
_{H^{s } (\mathbb{R} ^+ )} dt\right)^{\frac{1}{2}} \\ \ns
  \displaystyle \leq C   T^{ \frac{s}{4} -  \varepsilon +
\frac{1}{2}}       \left( \int_0^T  \| t^{ \frac{|s|}{4} +
\varepsilon} u_{xx}  \|^2 _{L^2 (\mathbb{R} ^+ )}
dt\right)^{\frac{1}{2}}
 \leq  C   T^{ \frac{s}{4} -  \varepsilon +
\frac{1}{2}}    \|u\| _{X_{s,T}^{\varepsilon}}.
\end{array}
\end{equation*}

{\it Step 2.} Now we prove that for $-2 + 4 \varepsilon < s<0$, it
holds
\begin{equation}\label{6.2.13}
\| u v _x \| _{L^1(0,T; H ^{s} (\mathbb{R} ^+ ))}    \leq C
T^{\alpha(s, \varepsilon)}
  \|u\|  _{X_{s,T}^{\varepsilon}} \|v\|  _{X_{s,T}^{\varepsilon}} .
\end{equation}
We distinguish four cases  $-1 \leq s <0$, $- \frac{3}{2} <s <-1$,
$s = - \frac{3}{2}$ and $-2 +4 \varepsilon <s < - \frac{3}{2} $.

 \medskip

{\it Case 1.}  $-1 \leq s <0$.  Noting the embedding $L^2
(\mathbb{R}^+) \subset H^{s} (\mathbb{R}^+)$  is continuous for
$s<0$, by  the Gagliardo-Nirenberg's inequality and the
H$\ddot{\text{o}}$lder's inequality, we get for any $u$, $v\in
X_{s,T}^{\varepsilon}$,
\begin{equation}\label{6.2.28}
\begin{array}{ll}
\displaystyle \| u v _x \| _{L^1(0,T; H ^{s} (\mathbb{R} ^+ ))}
\displaystyle \leq   C   \int_0^T  \left \|u   v _x  \right\| _{ L^2
(\mathbb{R} ^+ )} dt \\ \ns
\displaystyle  \leq   C  \int_0^T   t^{- \frac{|s|}{2} - 2
\varepsilon}    \|  t^{ \frac{|s|}{4} + \varepsilon }   u  \| _{ L^2
(\mathbb{R} ^+ )}    \|  t^{ \frac{|s|}{4} + \varepsilon }  v_x \|_{
L^{\infty} (\mathbb{R} ^+ )}   dt  \\ \ns
\displaystyle  \leq   C  \int_0^T   t^{- \frac{|s|}{2} - 2
\varepsilon}    \|  t^{ \frac{|s|}{4} + \varepsilon }   u  \| _{ L^2
(\mathbb{R} ^+ )}    \left(  \|  t^{ \frac{|s|}{4} + \varepsilon }
v_x \|_{ L^2 (\mathbb{R} ^+ )}      +      \| t^{\frac{|s|}{4} +
\varepsilon } v_x \| _{L^2( \mathbb{R}^+ )} ^{\frac{1}{2}}
\|t^{\frac{|s|}{4} + \varepsilon} v_{xx} \| _{L^2( \mathbb{R}^+ )}
  ^{\frac{1}{2}}     \right) dt  \\ \ns
 \displaystyle  =  C
\int_0^T   t^{- \frac{3|s|}{8} - \frac{3 \varepsilon}{2}}        \|
t^{ \frac{|s|}{4} + \varepsilon } u \| _{L^2( \mathbb{R}^+ )}
\| t^{\frac{|s|}{8} + \frac{\varepsilon}{2} } v_x \| _{L^2(
\mathbb{R}^+ )}  dt  \\ \ns
\displaystyle   \quad  +  C  \int_0^T
        t^{- \frac{7|s|}{16} - \frac{7 \varepsilon}{4}}       \| t^{ \frac{|s|}{4} + \varepsilon } u \|
_{L^2( \mathbb{R}^+ )}          \| t^{\frac{|s|}{8} +
\frac{\varepsilon}{2} } v_x \| _{L^2( \mathbb{R}^+ )} ^{\frac{1}{2}}
\|t^{\frac{|s|}{4} + \varepsilon} v_{xx} \| _{L^2( \mathbb{R}^+ )}
  ^{\frac{1}{2}}  dt\\ \ns
\displaystyle \leq  C   \sup_{t\in[0,T]} \| t^{\frac{|s|}{4} +
\varepsilon} u \| _{L^2( \mathbb{R}^+ )}      \Big[    \left(
\int_0^T t^{- \frac{3|s|}{4} - 3 \varepsilon}
dt\right)^{\frac{1}{2}} \left( \int_0^T     \|t^{ \frac{|s|}{8} +
\frac{\varepsilon}{2}} v _x \| _{L^2( \mathbb{R}^+ )}^2  dt  \right)
  ^{\frac{1}{2}} \\ \ns
\displaystyle \quad  +     \left( \int_0^T t^{- \frac{7|s|}{8} -
\frac{7 \varepsilon}{2} } dt\right)^{\frac{1}{2}} \left( \int_0^T
\|t^{ \frac{|s|}{8} + \frac{\varepsilon}{2}} v _x \| _{L^2(
\mathbb{R}^+ )}^2  dt \right)
  ^{\frac{1}{4}}      \left( \int_0^T     \|t^{
\frac{|s|}{4} + \varepsilon} v _{xx} \| _{L^2( \mathbb{R}^+ )}^2 dt
\right)
  ^{\frac{1}{4}} \Big].
\end{array}
\end{equation}
Again, by  the Gagliardo-Nirenberg's inequality and the
H$\ddot{\text{o}}$lder's inequality, one yields
\begin{equation}\label{6.2.31}
\begin{array}{ll}
\displaystyle    \int_0^T     \|t^{ \frac{|s|}{8} +
\frac{\varepsilon}{2}} v _x \| _{L^2( \mathbb{R}^+ )}^2  dt
=  \int_0^T    t^{ -\frac{|s|}{4} - \varepsilon}    \|t^{
\frac{|s|}{4} +  \varepsilon } v _x \| _{L^2( \mathbb{R}^+ )}^2  dt \\
\ns
\displaystyle    \leq C  \int_0^T     t^{ -\frac{|s|}{4} -
\varepsilon}    \left(  \| t^{ \frac{|s|}{4} +  \varepsilon } v \|_{
L^{2} (\mathbb{R} ^+ )} ^2    + \|t^{ \frac{|s|}{4} +  \varepsilon }
v \|_{ L^{2} (\mathbb{R} ^+ )}     \| t^{ \frac{|s|}{4} +
\varepsilon }  v_{xx} \|_{ L^{2} (\mathbb{R} ^+ )} \right) dt\\ \ns
\displaystyle \leq  C   \sup_{t\in[0,T]} \| t^{\frac{|s|}{4} +
\varepsilon} v \| _{L^2( \mathbb{R}^+ )} ^2        \int_0^T t^{-
\frac{|s|}{4} - \varepsilon } dt\\ \ns
\displaystyle \quad  +  C   \sup_{t\in[0,T]} \| t^{\frac{|s|}{4} +
\varepsilon} v \| _{L^2( \mathbb{R}^+ )}          \left( \int_0^T
t^{- \frac{|s|}{2} - 2\varepsilon } dt\right)^{\frac{1}{2}}
     \left( \int_0^T     \|t^{
\frac{|s|}{4} + \varepsilon} v _{xx} \| _{L^2( \mathbb{R}^+ )}^2 dt
\right)
  ^{\frac{1}{2}} .
\end{array}
\end{equation}
It follows from
\begin{equation*}
\left\{\begin{array}{ll}  \frac{3|s|}{4} + 3 \varepsilon < 1, \quad
\frac{7|s|}{8} + \frac{7 \varepsilon}{2}
< 1,  \\
\ns
  \frac{|s|}{4} + \varepsilon<1, \quad    \frac{|s|}{2} +
  2\varepsilon<1,
\end{array}\right.
\end{equation*}
that $ -\frac{8}{7} + 4 \varepsilon  < s <0$. Hence,  by
(\ref{6.1.1}), (\ref{6.2.28}) and (\ref{6.2.31}), we show  that
 estimate (\ref{6.2.13}) holds for $-1\leq s <0$ when $\varepsilon>0$ is small enough.

\medskip

{\it Case 2.} $-\frac{3}{2}  < s < -1$.  The operator
$\frac{\partial}{\partial x}:  H_0^{|s|} (\mathbb{R}^+)
 \rightarrow H^{|s| -1 } (\mathbb{R}^+) $ is linear bounded  for
$-\frac{3}{2}  < s < -1$. Then by transposition, we have
\begin{equation}\label{6.2.32}
 \int_0^T  \left \| (u  v)_x   \right\| _{ H ^{s} (\mathbb{R} ^+ )}  dt
\displaystyle \leq  C \int_0^T  \left \|u  v   \right\| _{ H ^{s+ 1}
(\mathbb{R} ^+ )} dt.
\end{equation}
Noting the embedding $L^{\frac{2}{2|s | -1} } (\mathbb{R}^+) \subset
H^{s + 1} (\mathbb{R}^+)$  is continuous for $-\frac{3}{2}  < s <
-1$, from (\ref{6.2.32}), by the H$\ddot{\text{o}}$lder's inequality
and the Gagliardo-Nirenberg's inequality, we get for any $u$, $v\in
X_{s,T}^{\varepsilon}$,
\begin{equation*}\label{6.2.33}
\begin{array}{ll}
\displaystyle \| ( u v) _x \| _{L^1(0,T; H ^{s} (\mathbb{R} ^+ ))}
\displaystyle
 \leq  C \int_0^T  \left \|u  v   \right\| _{ L^{\frac{2}{2|s | -1} } (\mathbb{R} ^+ )}
 dt \\ \ns
\displaystyle  = \int_0^T    t^{- \frac{|s|}{2} - 2 \varepsilon} \|(
t^{ \frac{|s|}{4} + \varepsilon}  u ) (  t^{ \frac{|s|}{4} +
\varepsilon} v  ) \| _{ L ^{\frac{2}{2|s| - 1} } (\mathbb{R} ^+ )}
dt
\\ \ns
\displaystyle \leq  C \int_0^T    t^{- \frac{|s|}{2} - 2
\varepsilon}    \| t^{ \frac{|s|}{4} + \varepsilon}  u \|_{L^2
(\mathbb{R} ^+ )}
   \|t^{ \frac{|s|}{4} + \varepsilon}  v\| _{ L ^{\frac{1}{|s| - 1} } (\mathbb{R} ^+ )}   dt
 \\ \ns
\displaystyle  \leq  C  \int_0^T   t^{- \frac{|s|}{2} - 2
\varepsilon} \| t^{ \frac{|s|}{4} +  \varepsilon } u \| _{L^2(
\mathbb{R}^+ )} \left( \| t^{\frac{|s|}{4} +\varepsilon } v \|
_{L^2( \mathbb{R}^+ )}          +                \| t^{\frac{|s|}{4}
+ \varepsilon } v \| _{L^2( \mathbb{R}^+ )} ^{\frac{|s |}{2} +
\frac{1}{4}} \|t^{\frac{|s|}{4} + \varepsilon} v_{xx} \| _{L^2(
\mathbb{R}^+ )}
  ^{-\frac{|s |}{2} + \frac{3}{4}}      \right)
  dt\\ \ns
\displaystyle \leq  C   \sup_{t\in[0,T]} \| t^{\frac{|s|}{4} +
\varepsilon} u \| _{L^2( \mathbb{R}^+ )}        \Big[
\sup_{t\in[0,T]} \| t^{\frac{|s|}{4} + \varepsilon  } v \| _{L^2(
\mathbb{R}^+ )}
     \int_0^T t^{- \frac{|s|}{2} - 2 \varepsilon}
dt  \\ \ns
\displaystyle \quad  +        \sup_{t\in[0,T]} \| t^{\frac{|s|}{4} +
\varepsilon  } v \| _{L^2( \mathbb{R}^+ )} ^{\frac{|s |}{2} +
\frac{1}{4}}         \left( \int_0^T  t^{- \frac{4|s|+ 16
\varepsilon }{2|s| +5} } dt\right)^{\frac{|s |}{4} + \frac{5}{8}}
\left( \int_0^T \|t^{ \frac{|s|}{4} + \varepsilon } v _{xx} \|
_{L^2( \mathbb{R}^+ )}^2 dt \right)
   ^{-\frac{|s |}{4} + \frac{3}{8}} \Big].
\end{array}
\end{equation*}
When $-\frac{3}{2}<s <-1$ and  $\varepsilon>0$ is small enough, it
holds
\begin{equation*}
\left\{\begin{array}{ll}  \frac{|s|}{2} + 2 \varepsilon < 1,
\\ \ns
\frac{4|s|+ 16 \varepsilon }{2|s| +5} < 1.
\end{array}\right.
\end{equation*}
Hence, combining (\ref{6.1.1}), one concludes estimate
(\ref{6.2.13}) for $-\frac{3}{2}  < s < -1$.

\medskip

{\it Case 3.} $s= -\frac{3}{2}$. Noting the embedding $H^{-1}
(\mathbb{R}^+) \subset H^{-\frac{3}{2} } (\mathbb{R}^+)$  is
continuous, by  the Gagliardo-Nirenberg's inequality and the
H$\ddot{\text{o}}$lder's inequality,
\begin{equation}\label{6.2.29}
\begin{array}{ll}
\displaystyle \| ( u v )_x \| _{L^1(0,T; H ^{-\frac{3}{2}}
(\mathbb{R} ^+ ))}
\displaystyle \leq   C   \int_0^T  \left \| (u   v) _x  \right\| _{
H^{-1} (\mathbb{R} ^+ )} dt \\  \ns
\displaystyle  =   C   \int_0^T    t^{- \frac{|s|}{2} - 2
\varepsilon }   \|  ( t^{ \frac{|s|}{2} + 2 \varepsilon }  u v) _x
\| _{ H^{-1} (\mathbb{R} ^+ )} dt
\displaystyle  \leq  C  \int_0^T  t^{ -\frac{|s|}{2} - 2 \varepsilon
}  \|  t^{ \frac{|s|}{2} + 2 \varepsilon }  uv  \| _{ L^2
(\mathbb{R} ^+ )} dt
\\ \ns
\displaystyle \leq   C  \int_0^T    t^{- \frac{|s|}{2} - 2
\varepsilon }    \|   t^{ \frac{|s|}{4} + \varepsilon }   u  \| _{
L^2 (\mathbb{R} ^+ )}     \|t^{ \frac{|s|}{4} + \varepsilon } v \|_{ L^{\infty} (\mathbb{R} ^+ )}   dt  \\
\ns
\displaystyle    \leq C  \int_0^T   t^{- \frac{|s|}{2} - 2
\varepsilon }    \| t^{\frac{|s|}{4} + \varepsilon } u  \| _{ L^2
(\mathbb{R} ^+ )}     \left(  \| t^{\frac{|s|}{4} + \varepsilon } v
\|_{ L^{2} (\mathbb{R} ^+ )} + \| t^{\frac{|s|}{4} + \varepsilon } v
\|_{ L^{2} (\mathbb{R} ^+ )} ^{\frac{1}{2}}   \| t^{\frac{|s|}{4} +
\varepsilon } v_x \|_{ L^{2} (\mathbb{R}
^+ )} ^{\frac{1}{2}}    \right)  dt  \\
\ns
\displaystyle  =  C  \int_0^T   t^{- \frac{|s|}{2} - 2 \varepsilon}
\| t^{ \frac{|s|}{4} + \varepsilon } u \| _{L^2( \mathbb{R}^+ )}
\| t^{\frac{|s|}{4} + \varepsilon } v \| _{L^2( \mathbb{R}^+ )}  dt \\
\ns
\displaystyle  \quad  + C  \int_0^T   t^{- \frac{7|s|}{16} - \frac{7
\varepsilon}{4}}    \| t^{ \frac{|s|}{4} + \varepsilon } u \| _{L^2(
\mathbb{R}^+ )}    \| t^{ \frac{|s|}{4} + \varepsilon } v
\|^{\frac{1}{2}}  _{L^2( \mathbb{R}^+ )} \|t^{\frac{|s|}{8} +
\frac{\varepsilon}{2} } v_{x} \| _{L^2( \mathbb{R}^+ )}
  ^{\frac{1}{2}} dt\\ \ns
\displaystyle \leq  C   \sup_{t\in[0,T]} \| t^{\frac{|s|}{4} +
\varepsilon} u \| _{L^2( \mathbb{R}^+ )}     \Big[  \sup_{t\in[0,T]}
\| t^{\frac{|s|}{4} + \varepsilon} v \| _{L^2( \mathbb{R}^+ )}
\int_0^T t^{- \frac{|s|}{2} - 2 \varepsilon} dt \\ \ns
\displaystyle \quad  +        \sup_{t\in[0,T]} \| t^{\frac{|s|}{4} +
\varepsilon} v \| _{L^2( \mathbb{R}^+ )} ^{\frac{1}{2}}
  \left( \int_0^T t^{- \frac{7|s|}{12} - \frac{7 \varepsilon}{3} }
dt\right)^{\frac{3}{4}}          \left( \int_0^T \|t^{ \frac{|s|}{8}
+ \frac{\varepsilon}{2}} v _x \| _{L^2( \mathbb{R}^+ )}^2  dt
\right)
  ^{\frac{1}{4}}    \Big] .
\end{array}
\end{equation}
When $s= - \frac{3}{2}$ and  $\varepsilon>0$ is small enough, it
holds
\begin{equation*}
\left\{\begin{array}{ll}   \frac{|s|}{2} +
  2\varepsilon<1, \quad
\frac{7|s|}{12} + \frac{7 \varepsilon}{3}
< 1,  \\
\ns
  \frac{|s|}{4} + \varepsilon<1.
\end{array}\right.
\end{equation*}
Hence, by (\ref{6.1.1}),  (\ref{6.2.31}) and (\ref{6.2.29}), one has
that estimate (\ref{6.2.13}) holds for $s= - \frac{3}{2}$.

\medskip

{\it Case 4.}  $-2 + 4 \varepsilon < s < -\frac{3}{2}$. The operator
$\frac{\partial}{\partial x}:  H_0^{|s|} (\mathbb{R}^+)
 \rightarrow H_0^{|s| -1 } (\mathbb{R}^+) $ is linear bounded  for
$-2 + 4 \varepsilon  < s < -\frac{3}{2}$. Then by transposition, we
have
\begin{equation}\label{6.2.39}
 \int_0^T  \left \| (u  v)_x   \right\| _{ H ^{s} (\mathbb{R} ^+ )}  dt
\displaystyle \leq  C \int_0^T  \left \|u  v   \right\| _{ H ^{s+ 1}
(\mathbb{R} ^+ )} dt.
\end{equation}
Noting the embedding $L^{1} (\mathbb{R}^+) \subset H^{s + 1}
(\mathbb{R}^+)$  is continuous, from (\ref{6.2.39}), by the
H$\ddot{\text{o}}$lder's inequality and the Gagliardo-Nirenberg's
inequality, we get for any $u$, $v\in X_{s,T}^{\varepsilon}$,
\begin{equation*}\label{6.2.36}
\begin{array}{ll}
\displaystyle \| ( u v) _x \| _{L^1(0,T; H ^{s} (\mathbb{R} ^+ ))}
\leq  C \int_0^T  \|u  v  \| _{ L^1 (\mathbb{R} ^+ )} dt\\
\ns
\displaystyle =  C \int_0^T   t^{- \frac{|s|}{2} - 2 \varepsilon}
  \| ( t^{ \frac{|s|}{4} +  \varepsilon} u )   (t^{ \frac{|s|}{4} +  \varepsilon}   v ) \| _{ L^1 (\mathbb{R} ^+ )}   dt
 \\ \ns
\displaystyle  \leq  C  \int_0^T   t^{- \frac{|s|}{2} - 2
\varepsilon} \| t^{ \frac{|s|}{4} +  \varepsilon } u \| _{L^2(
\mathbb{R}^+ )}
     \| t^{\frac{|s|}{4} +\varepsilon } v \| _{L^2( \mathbb{R}^+
)}  dt\\ \ns
\displaystyle \leq  C   \sup_{t\in[0,T]} \| t^{\frac{|s|}{4} +
\varepsilon} u \| _{L^2( \mathbb{R}^+ )}      \sup_{t\in[0,T]} \|
t^{\frac{|s|}{4} + \varepsilon  } v \| _{L^2( \mathbb{R}^+ )}
     \int_0^T t^{- \frac{|s|}{2} - 2 \varepsilon}
dt  .
\end{array}
\end{equation*}
It follows from $ \frac{|s|}{2} + 2 \varepsilon < 1$ that  $ -2 + 4
\varepsilon < s $.
Hence, estimate (\ref{6.2.13}) holds for $ -2 + 4 \varepsilon < s <
-\frac{3}{2}$.

\medskip

{\it Step 3.} Claim for any $u$, $v \in X_{s,T}^{\varepsilon}$,
\begin{equation}\label{6.2.25}
 \| t^{\frac{|s|}{4} + \varepsilon}  u v_x   \| _{L^1(0,T; L^2 (\mathbb{R} ^+ ))} \leq C
T^{\alpha(s, \varepsilon)}
   \|u\|  _{X_{s,T}^{\varepsilon}} \|v\|  _{X_{s,T}^{\varepsilon}} .
\end{equation}
In fact,  from the Gagliardo-Nirenberg's inequality and the
H$\ddot{\text{o}}$lder's inequality, we obtain that for any $u$,
$v\in X_{s,T}^{\varepsilon}$,
\begin{equation*}\label{6.2.27}
\begin{array}{ll}
\displaystyle  \| t^{\frac{|s|}{4} + \varepsilon}  u v_x   \|
_{L^1(0,T; L^2 (\mathbb{R} ^+ ))}
\displaystyle =  \int_0^T  t^{- \frac{|s|}{8} -
\frac{\varepsilon}{2}}   \| ( t^{ \frac{|s|}{4} + \varepsilon} u )
 ( t^{ \frac{|s|}{8} + \frac{\varepsilon}{2}} v )_x \| _{L^2(\mathbb{R} ^+ )} dt\\ \ns
\displaystyle  \leq  \int_0^T  t^{- \frac{|s|}{8} -
\frac{\varepsilon}{2}}    \| t^{ \frac{|s|}{4} + \varepsilon} u \|
_{L^{\infty} (\mathbb{R} ^+ )}
 \| t^{ \frac{|s|}{8} + \frac{\varepsilon}{2}} v _x  \| _{L^2 (\mathbb{R} ^+ )} dt \\
\ns
\displaystyle  \leq  C  \int_0^T  t^{- \frac{|s|}{8} -
\frac{\varepsilon}{2}} \| t^{\frac{|s|}{4} + \varepsilon} u \|
_{L^2( \mathbb{R}^+ )}     \| t^{ \frac{|s|}{8} +
\frac{\varepsilon}{2}} v _x \| _{L^2( \mathbb{R}^+ )}
  dt\\ \ns
\displaystyle  \quad  +  C  \int_0^T  t^{- \frac{|s|}{8} -
\frac{\varepsilon}{2}} \| t^{\frac{|s|}{4} + \varepsilon} u \|
_{L^2( \mathbb{R}^+ )} ^{\frac{1}{2}} \|t^{\frac{|s|}{4} +
\varepsilon} u_x \| _{L^2( \mathbb{R}^+ )}
  ^{\frac{1}{2}}   \| t^{ \frac{|s|}{8} + \frac{\varepsilon}{2}} v _x \| _{L^2( \mathbb{R}^+ )}
  dt\\ \ns
\displaystyle \leq  C   \sup_{t\in[0,T]} \| t^{\frac{|s|}{4} +
\varepsilon} u \| _{L^2( \mathbb{R}^+ )}          \left( \int_0^T
t^{- \frac{|s|}{4} - \varepsilon} dt\right)^{\frac{1}{2}}    \left(
\int_0^T     \|t^{ \frac{|s|}{8} + \frac{\varepsilon}{2}} v _x \|
_{L^2( \mathbb{R}^+ )}^2  dt  \right)
  ^{\frac{1}{2}} \\ \ns
\displaystyle \quad  +  C   \sup_{t\in[0,T]} \| t^{\frac{|s|}{4} +
\varepsilon} u \| _{L^2( \mathbb{R}^+ )}  ^{\frac{1}{2}} \\ \ns
\displaystyle  \quad  \left( \int_0^T t^{- \frac{|s|}{2} -
2\varepsilon} dt\right)^{\frac{1}{4}} \left( \int_0^T
\|t^{\frac{|s|}{4}+ \varepsilon} u_x \| _{L^2( \mathbb{R}^+ )}^2 dt
 \right)^{\frac{1}{4}} \left( \int_0^T     \|t^{ \frac{|s|}{8} + \frac{\varepsilon}{2}}
v _x \| _{L^2( \mathbb{R}^+ )}^2  dt  \right)
  ^{\frac{1}{2}} ,
\end{array}
\end{equation*}
which combining  (\ref{6.1.1}) and (\ref{6.2.31})  implies estimate
(\ref{6.2.25}).
This completes the proof of Lemma \ref{Lma 5.1}.\endpf

\begin{proposition}\label{Pro 6.5}
Let  $0< T\leq 1$,  $\varepsilon>0$,  $-2 + 4 \varepsilon < s< 0$,
$\phi \in H^s (\mathbb{R}^+)$, $\vec{h} \in {\cal H} ^s ( 0,T )$ and
$t^{ \frac{|s|}{4} + \varepsilon} \vec{h} \in {\cal H} ^{0} (0,T)$.
Then there exists a $T_* \in (0,T]$ depending on $ \| \phi \|
_{H^s(\mathbb{R}^+)}   + \| \vec{h} \| _{{\cal H} ^s (0,T)} +   \|
t^{\frac{|s|} {4} + \varepsilon } \vec{h} \| _{ {\cal H}^0 (0,T) }
$, such that equation (\ref{1.1.1*}) admits a unique solution  $u
\in X_{s,T_*}^{\varepsilon} $.
Moreover, the corresponding solution map from the space of initial
and boundary data to the solution space is  continuous.

\end{proposition}

{\it Proof. } The solution of system (\ref{1.1.1*}) can be written
in the form
\begin{equation*}\label{6.3.3}
u(t) = W_c(t) \phi + W_{bdr} (t) \vec{h}  - \int_0^t W_c(t-\tau)
\left( u_{xx} + u u_x  \right) (\tau)d\tau .
\end{equation*}
For $w\in X_{s,T_*}^{\varepsilon} (d)  : = \left \{w\in
X_{s,T_*}^{\varepsilon}   \; \big|\ \|w\| _{X_{s,T_*}^{\varepsilon}
} \leq d \right \}$, define
\begin{equation}\label{6.3.7}
\Gamma (w) = W_c(t) \phi + W_{bdr} (t) \vec{h}  - \int_0^t
W_c(t-\tau) \left( w_{xx} + w w_x  \right) (\tau)d\tau.
\end{equation}
Then by Proposition \ref{Pro 2.1*} - Proposition \ref{Pro 2.3*},
Lemma \ref{Lma 5.1} and (\ref{6.3.7}), we obtain that for any $w$,
$w_1$, $w_2 \in X_{s,T_*}^{\varepsilon} (d) $,
\begin{equation*}\label{6.3.19}
\begin{array}{ll}
\displaystyle
  \|\Gamma (w) \|_{ X_{s,T_*}^{\varepsilon} } \\  \ns
\displaystyle \leq   C_1  \left( \|  \phi \| _{H^s(\mathbb{R} ^+ )}
+  \| \vec{h} \|_{ {\cal H}^s (0,T) } + \| t^{\frac{|s|} {4} +
\varepsilon } \vec{h} \| _{ {\cal H}^0 (0,T) } \right)
  +     C_2  T_*^{\alpha(s,\varepsilon)} \left(
  \|w\|  _{X_{s,T_*}^{\varepsilon}}   +
 \|w\| _{X_{s,T_*}^{\varepsilon}}  ^2 \right),
\end{array}
\end{equation*}
and
\begin{equation*}\label{6.3.10}
\begin{array}{ll}
\|  \Gamma(w_1)   -   \Gamma(w_2)  \|_{X_{s,T_*}^{\varepsilon}} \\
\ns
\displaystyle  \leq    C_2    \left(  T_*^{\alpha (s, \varepsilon)}
  \|w_1 - w_2\|  _{X_{s,T_*}^{\varepsilon}}   +   T_*^{\alpha (s, \varepsilon)}
 \|w_1 + w_2\|_{X_{s,T_*}^{\varepsilon}}  \|w_1 - w_2\|_{X_{s,T_*}^{\varepsilon}}  \right) .
 \end{array}
\end{equation*}
Similar to the proof of Proposition \ref{Pro 3.1}, equation
(\ref{1.1.1*}) admits a unique solution $u = \Gamma(u)$ in $
X_{s,T_*}^{\varepsilon} (d) $.\endpf

\smallskip

\medskip

\noindent {\bf Proof of Theorem \ref{Thm 1.2*} for $-2 < s<0$. }
Proposition \ref{Pro 6.5} shows that there exists a $T_*\in (0, 1]$,
such that equation (\ref{1.1.1*}) has a local solution $u\in
C([0,T_*]; H^s (\mathbb{R}^+) )$  with $t^{ \frac{|s|}{4} +
\varepsilon} u  \in C([0,T_*]; L^2 (\mathbb{R}^+) ) $.   Then
$u(x,T_*) \in L^2 (\mathbb{R}^+)$. From $t^{ \frac{|s|}{4} +
\varepsilon} \vec{h} \in {\cal H} ^{0} (0,T)$, we see that $\vec{h}
\in {\cal H} ^{0} (T_*,T)$.
By the result of Theorem \ref{Thm 1.2*} for $s=0$, equation
(\ref{1.1.1*}) with initial datum $u(x,T_*) \in L^2 (\mathbb{R}^+)$
and boundary conditions $\vec{h} \in {\cal H} ^{0} (T_*,T)$ admits a
unique solution $u\in C([T_*, T]; H^s (\mathbb{R}^+) )$.
\endpf



\begin{thebibliography}{99}




\bibitem{Aimar} M.~T.~Aimar, \sl Etude num$\acute{e}$rique d'une $\acute{e}$quation d$\acute{e}$volution non lin$\acute{e}$aire d$\acute{e}$rivant l'instabilit$\acute{e}$ thermodiffusive d'un front de flamme,
 \rm Th$\grave{e}$se 3 $\acute{e}$me cycle Universit$\acute{e}$ de Provence,
1982.

\bibitem{Armaou} A.~Armaou and P.~D.~Christofides, \sl Feedback control of the Kuramoto-Sivashinsky
equation, \sl Phys. D,
 \rm{\bf 137} (2000), 49--61.


\bibitem{Biagioni} H.~A.~Biagioni and T.~Gramchev, \it Multidimensional Kuramoto-Sivashinsky type equations:
singular initial data and analytic regularity,
 \sl Fifth Workshop on Partial Differential Equations (Rio de Janeiro,
1997). Mat. Contemp.,  \rm{\bf 15} (1998), 21--42.



\bibitem{Bona} J.~L.~Bona and L.~R.~Scott, \it Solutions of the Korteweg-de Vries equation in fractional order Sobolev spaces,
 \sl Duke  Math. J.,  \rm{\bf 43}
(1976), 87--99.

\bibitem{Zhang1} J.~L.~Bona, S.~Sun and B.-Y.~Zhang, \it A non-homogeneous boundary-value problem for the Korteweg-de Vries equation in a quarter plane,
 \sl Trans. Amer. Math. Soc.,  \rm{\bf 354}
(2002), 427--490.


\bibitem{Zhang2} J.~L.~Bona, S.~Sun and B.-Y.~Zhang, \it A nonhomogeneous boundary-value problem for the Korteweg-de Vries equation posed on a finite domain,
 \sl Comm. Partial Differential Equations,  \rm{\bf 28}
(2003), 1391--1436.

\bibitem{Zhang3} J.~L.~Bona, S.~Sun and B.-Y.~Zhang, \it A non-homogeneous boundary-value problem for the Korteweg-de Vries equation posed on a finite domain II,
 \sl J. Differential Equations,  \rm{\bf 247}
(2009), 2558--2596.



\bibitem{Cerpa1}
E.~Cerpa, \it Null controllability and stabilization of the linear
Kuramoto-Sivashinsky equation, \sl Commun. Pure Appl. Anal., \rm{\bf
9} (2010), 91--102.


\bibitem{Cerpa2}
E.~Cerpa and A.~Mercado, \it Local exact controllability to the
trajectories of the 1-D Kuramoto-Sivashinsky equation, \sl J.
Differential Equations, \rm{\bf 250} (2011), 2024--2044.

\bibitem{Collet}
P.~Collet, J.-P.~Eckmann, H.~Epstein and J.~Stubbe, \it A global
attracting set for the Kuramoto-Sivashinsky equation, \sl Comm.
Math. Phys., \rm{\bf 152} (1993), 203--214.

\bibitem{Cousin}
A.~T.~Cousin and N.~A.~Larkin, \it Kuramoto-Sivashinsky equation in
domains with moving boundaries, \sl Port. Math., \rm{\bf 59} (2002),
335--349.

\bibitem{Demirkaya}
A.~Demirkaya, \it The existence of a global attractor for a
Kuramoto-Sivashinsky type equation in 2D, \sl  Discrete Contin. Dyn.
Syst., \rm (2009), 198--207.


\bibitem{Giacomelli}
L.~Giacomelli and F. Otto, \it New bounds for the
Kuramoto-Sivashinsky equation, \sl  Commun. Pure Appl. Math.,
\rm{\bf 58}  (2005), 297--318.


\bibitem{Glowinski} R.~Glowinski, J.~L.~Lions and J. He, \sl Exact and Approximate Controllability for Distributed Parameter
Systems: A Numerical Approach, \rm Encyclopedia of Mathematics and
its Applications,  117, Cambridge University Press, Cambridge, 2008.




\bibitem{Ilyashenko}
Ju.~S.~Il'yashenko, \it Global Analysis of the Phase Portrait for
the Kuramoto-Sivashinsky Equation, \sl J. Dynam. Differential
Equations, \rm{\bf 4} (1992), 585--615.



\bibitem{Kaikina} E.~I.~Kaikina, \it Subcritical Kuramoto-Sivashinsk-type equation on a half-line,
\sl J.  Differential Equations,  \rm {\bf 220} (2006), 279--321.




\bibitem{Kuramoto2} Y.~Kuramoto, \it On the formation of dissipative structures in reaction-diffusion systems, \sl
Progr. Theoret. Phys.,  \rm {\bf 54} (1975), 687--699.

\bibitem{Kuramoto3} Y.~Kuramoto, \it Instability and turbulence of wavefronts in reaction-diffusion systems, \sl
Progr. Theoret. Phys.,  \rm {\bf 63} (1980), 1885--1903.


\bibitem{Lions} J.~L.~Lions and E. Magenes, \sl Non-Homogeneous Boundary Value Problems and Applications, \rm vol. II,
Springer-Verlag,  Berlin-Heidelberg-New York, 1972.
%


\bibitem{Nicolaenko1} B.~Nicolaenko and B.~Scheurer, \it Remarks on the Kuramoto-Sivashinsky equation, \sl Phys. D,
 \rm{\bf 12} (1984), 391--395.

\bibitem{Nicolaenko2} B.~Nicolaenko, B.~Scheurer and R.~Temam, \it Some global dynamical properties of the Kuramoto-Sivashinsky equations: nonlinear stability and attractors,
\sl Phys. D, \rm{\bf 16} (1985), 155--183.

\bibitem{Pazy}  A.~Pazy, \sl Semigroups of Linear Operators and Applications to Partial Differential Equations,
\rm Applied Mathematical Sciences, vol. 44,
Springer-Verlag, New York, 1983.

\bibitem{Pilod}  D.~Pilod,
\it Sharp well-posedness results for the Kuramoto-Velarde equation,
\sl Commun. Pure Appl. Anal., \rm{\bf 7} (2008), 867--881.


\bibitem{Sivashinsky1}  G.~I.~Sivashinsky,
\it Nonlinear analysis of hydrodynamic instability in laminar flames--I. Derivation of basic equations, \sl Acta Astronaut., \rm{\bf 4} (1977),
1177--1206.

\bibitem{Sivashinsky2}  G.~I.~Sivashinsky, \it On flame propagation under conditions of stoichiometry, \sl SIAM J. Appl. Math.,  \rm{\bf 39}
(1980), 67--82.


\bibitem{Tadmor} E.~Tadmor, \it The well-posedness of the Kuramoto-Sivashinsky equation, \sl SIAM J. Math. Anal.,  \rm{\bf 17}
(1986), 884--893.


\bibitem{Tartar} L.~Tartar, \it Interpolation non lin$\acute{e}$aire et r$\acute{e}$gularit$\acute{e}$,
\sl J. Funct. Anal.,  \rm{\bf 9} (1972), 469--489.




\end{thebibliography}
\end{document}